

\input amstex
\documentstyle{amsppt}
\magnification =1100
\hsize =14truecm
\hcorrection{0.6in}
\vcorrection{0.2in}
\NoBlackBoxes
\catcode`\@=11
\outer\def\Refs{\relaxnext@\def\refskip@{\hskip\@ne sp\hskip\m@ne sp}%
\def\next@{\ifx\next\nofrills\def\next@\nofrills{\nextii@}\else
\def\next@{\nextii@{References}}\fi\next@}%
\def\nextii@##1{\bigbreak\hbox to\hsize{\hfil\tenpoint
\smc\ignorespaces##1\unskip\hfil}\nobreak
\bigskip\tenpoint\sfcode`.=\@m}
\futurelet\next\next@}
\catcode`\@=\active%
\vsize =20.6truecm
\rightheadtext\nofrills{The singularities of the parameter surface}
\leftheadtext\nofrills{A. Grassi}
\topmatter
\title
The singularities of the parameter surface of a minimal elliptic threefold
\endtitle
\author
Antonella Grassi
\endauthor
\affil
Tufts University
\endaffil
\address
Department of Mathematics,  Tufts University, Medford, MA 02155
\endaddress
\email
agrassi\@jade.tufts.edu
\endemail
\thanks
 I would like to thank J. Koll\'ar and D. Morrison for offering valuable
suggestions\endthanks
\endtopmatter

\def\bda{\bold\Delta}
\def\bla{\bold\Lambda}
\def\bsa{\bold\Sigma}

\document

 Let $\pi : X \rightarrow S $  be an elliptic fibration
between smooth varieties.
If $\dim(X)=2$, then $X$ can be assumed minimal and  $K_X$ can be expressed as
the pullback of a divisor
from $S$.  Moreover $S$ cannot be changed birationally.    These observations
were the first step in Kodaira's analysis of
singular fibers [Kd].
In higher dimensions the birational geometry of the total space is more
complicated and the parameter space can be modified birationally.

 If $X $ has dimension $3$ and is not uniruled, then $X$ has a minimal model
(with terminal singularities) [Mo].
 In earlier work [G1] we have shown that there exists a birationally equivalent
elliptic fibration
$\bar {\pi} : \bar X  \rightarrow  \bar S $ such that  $\bar X$ is minimal and
$K_{\bar X} \equiv \bar {\pi} ^* ( K_{\bar S} + \bar {\bla}) $ (\S 0), [G1].
Moreover $\bar S$ has at worst quotient singularities;
 it is not difficult to find examples where $\bar S$ is actually singular.

\smallpagebreak

In this paper we study the singularities of this parameter surface $\bar S$.
Although   $\bar S$ is not uniquely determined
by the birational equivalence class of the fibration, any two
 such $\bar S$'s are related by  a particular kind of birational map (0.4).

The following result covers the case of no multiple fibers.

\proclaim{ Theorem 4.5} Let $\pi  : X \rightarrow S$ and $\bar S$ be as above
and assume that $\pi$ does not have multiple fibers outside a set of
codimension 2.  Then the only exceptional graphs of the minimal resolution of
$\bar S$ are Hirzebruch-Jung strings of type ${<} n , q {>}, \ n \leq 6$.

 In particular all the possible dual graphs are represented by the full circles
in Tables 2, 3, 4.
\endproclaim

\noindent In \S 3 we present examples of  all these ${<} n , q {>}$ other than
$n = 5, \ q = 2,3, 4$.

\smallpagebreak

If we allow multiple fibers things are more complicated.

 Kawamata [Ka5] has shown that any quotient singularity can occur on $\bar S$.
In his examples $\bar S$ has only one singular point $p$, $\bar \pi ^{-1} (p)$
is a multiple fiber, and the elliptic fibration is isotrivial outside $p$.
Moreover Dolgachev has shown me an example where $\bar S$ has a singular point
$p$, $\bar \pi ^{-1} (p)$ is a multiple fiber, and the elliptic fibration is
not isotrivial outside $p$. However also in this case the only divisors of
singular fibres passing through $p$ are the ones of multiple fibers.

\smallpagebreak

The following corollary is a consequence of Theorem 4.5.

We start with any elliptic fibration  $ X _0 \rightarrow S_0 $ and take a
suitable birational equivalent fibration $X \rightarrow S$
of smooth varieties (see 0.3).
Then  the birational map  $\phi: S  \rightarrow \bar S $ is actually a
morphism.

\proclaim{ Corollary 4.7} Assume that all the exceptional curves for $\phi$
which are image of multiple fiber divisors have zero intersection with  curves
image of  non multiple singular divisor.
 Then:
$$ \chi ( \Cal O _X) = - \frac{1}{2} \{ [ K _ {\bar S  } + \phi _* \pi _*
(K_{X/S}) ]
      \cdot \phi_* \pi (K_{X/S}) - \sum \Sb p \in V \endSb (1 - 1/ {n_p}) \},
\tag 4.7$$
 where $ V = \{ p \in U \text{ such that the fiber over the general point of  }
\phi ^{-1} (p) \text{ is not multiple } \}$.

\endproclaim
$\phi_* \pi _* (K_{X/S}) $ does not depend on the original choice of the
birational model.

This formula generalizes to the 3-dimensional case  the formula for the Euler
characteristic for an elliptic surface $p: S \rightarrow C$, namely:  $ \chi
(\Cal O _S)  =  \deg \  p_* K_{S/C}$.
If the fibration does not have multiple fiber, then
$\chi ( \Cal O _X) = - \frac{1}{2} \{( K _ {\bar S  } + \bar\bla ) \cdot
\bar\bla - \sum \Sb p \in V \endSb (1 - 1/ {n_p}) \} $.

  For example if  the Kodaira dimension of $X$ is equal to $0$, then $\chi (
\Cal O _X) =
 1 / 2   \sum \Sb p \in V \endSb (1 - 1/ {n_p})$.
   ($K_{\bar X} \equiv \bar {\pi} ^* ( K_{\bar S} + \bar {\bla})$, see \S 0.)
Other applications of this results can be found in [G1]

\smallpagebreak

The hypothesis of the above Corollary is  always satisfied if $\pi$ does not
have multiple fibers outside a set of codimension 2.

\vskip 0.2in

\heading {\bf 0.  Background results, notations and first properties. }
\endheading

\medpagebreak

\subheading{Background}

\smallpagebreak

Let $\pi :   X \rightarrow S$ be an elliptic fibration between smooth
varieties.

\definition{Definition} Set:

\noindent $\bsa _{X/S} = \{  P \in S  \ \text{ such that } \ \pi \text{ is not
smooth over } P \}$

\noindent \phantom{$\bsa _{X/S} = $}(the ramification divisor of $\pi$)

\noindent $\bsa ^ {om} _{X/S} = \{  P \in \bsa _{X/S}  \ \text{ such that } \
\pi  ^{-1} (P) \text{ is a multiple fiber} \} $ and  $\bsa ^ m _{X/S}$ the
closure of  $\bsa ^ {om} _{X/S} $.

We write $\bsa$ ({\rm resp. } $\bsa ^m$) for $\bsa _{X/S}$ ({\rm resp  } $\bsa
^ m _{X/S}$) when the fibration is well understood.
\enddefinition

We denote by $ J_\infty$ the divisor of poles of the elliptic modular function
$\bold J$.

The following theorems hold:

\proclaim{\bf Theorem 0.1  [F2]} Let  $\pi : X \rightarrow S$ be an elliptic
fibration between smooth varieties.
Assuming that the modular function $\bold J$ extends to an holomorphic map $S
\rightarrow \Bbb P^1$, then
$$ mK_ X  = \pi ^*( m K_S + m \pi_* (K_{X  /S}) + m\sum
( \frac{m_i - 1}{ m_i}  Y_i) )+ mE -  mG $$
where the fiber over the general point of
$Y_i$ is a multiple fiber of multiplicity $m_i$, $m$ is a multiple of the $ \{
m_i \}$'s; $mE$  and $mG$ are effective divisors.
\endproclaim

\proclaim{\bf Theorem 0.2 [Ka2]} Let  $\pi :   X \rightarrow S$ be an elliptic
fibration between smooth varieties.  Suppose that the ramification locus is a
divisor with simple normal crossings.
Then the modular function $\bold J$ extend to a holomorphic map
$S \rightarrow \Bbb P^1$, $12\pi _* (K_{X/S}) $ is an invertible sheaf and $
12\pi _* (K_{X/S}) \cong \Cal O _S (\sum 12 a_i D_i ) \otimes J _\infty.$

\noindent $D_i$ are irreducible components of $\bsa$ and $0 \leq a_i < 1$ the
rational numbers corresponding to the type
of singularities over the general point of $D_i$ and are described explicitly
in [Kd].

 If we write  $ J _ \infty = \sum b_j  B _j $ , then $\bold J$ has a pole of
order $ b _j    $ along $B _j $.
\endproclaim

For convenience, we give the following:

\definition{Definition}
Let $\pi :   X \rightarrow S$ be an elliptic fibration between smooth
varieties. Assume that $\bsa$ is a divisor with simple normal crossings.
Define:
$$ \aligned \bda  _{X/S} & =  \sum 12 a_i D_i + (1/12) J _\infty\\
\bla _{X/S} &=  \bda  _{X/S} + \sum \frac{m_i - 1}{ m_i}  Y_i \text{; }
\endaligned
 $$
When there is no danger of confusion, we denote these simply by $\bda$ and
$\bla$ respectively.  $\bla$ is thus a $\Bbb Q$-divisor.

\enddefinition
Note that $[ \bla] = 0$, where $[x]$ denotes the integral part of $x$ (see
0.2);
$12 m \bla _{X/S}$ is a divisor supported on $\bsa _{X/S}$.

In this notation we have

$$ \phantom{mnb}mK_ X  = \pi ^* m( K_S + \bla ) + mE -  mG. \tag  0.1'$$

\medpagebreak

\proclaim{ \bf Theorem  0.3 [G1]}
Let $X _0 \rightarrow S_0$ be an elliptic threefold which is not uniruled.

 Then there exist a birationally equivalent fibration
$\bar \pi : \bar X \rightarrow \bar S$, such that $\bar X$ has at worst
terminal and $\bar S$ log terminal singularities. Furthermore $K_{\bar X }$ is
nef and $K_{\bar X } \equiv \bar \pi ^* (K_ {\bar S}
 +  \bar\bla)$; $\bar\bla$ is the $\Bbb Q$-divisor defined below.
\endproclaim

\demo{\it Sketch of the proof} Let $\tilde X \rightarrow \tilde S$ be    a
birational equivalent fibration of smooth varieties.
 Let $ X _1 \rightarrow S $ be a flat model of the fibration  and $ X $ a non
singular model of $X_1$. Denote by $\pi : X \rightarrow S$ the induced
fibration.
Assume that  $\bsa _{X/S}$ (the ramification divisor of $\pi$) is a divisor
with simple normal crossings.

Since $X$ is not uniruled, $\kappa (X) $ is non-negative [Mi] and on some
suitable birational model, some multiple of $ K_S + \bla $ is an effective
divisor.
  Then there exists a unique surface $\bar S$, a divisor $\bar\bla$ and a
morphism $ \phi  : S \rightarrow \bar S$ such that $ K_{\bar S} + \bar\bla$ is
nef and $\phi ^* (K_{\bar S} + \bar\bla) + \sum c_i  \Gamma _i$ is the Zariski
decomposition of $K_S + \bla$ [Za], [Sa].  Furthermore $\bar\bla = \phi _* (
\bla)$,  and $\bar S$ has at worst log terminal singularities:
 $\phi$ is a {\it log terminal} contraction, because $[\bla] =0$ [KMM, \S
5.1.5.].

In particular  if $K_S + \bla$ is not nef, then there exists a {\it log
extremal}  curve  $\Gamma _1$ such that $(K_S + \bla ) \cdot \Gamma _1 <0$, and
a morphism
 $\phi _1 : S \to S_1$.  We successively construct the objects $\phi _i  :
S_{i-1}  \to S_i$ where $\phi _i$ is the log extremal contraction associated to
the log extremal curve $\Gamma _i$. All such $\Gamma _i$ are  rational curves
[Ka4], [Br].
$\phi$ is the composition of these log extremal contractions (a log
contraction).
Now we can apply relative Mori's theory [Mo] to the fibration $X \rightarrow
\bar S$, and this concludes the proof.
\qed
\enddemo

Given $S$ and $\bla$ as in (0.3), $\bar S$ is uniquely determined [Sa, Thm
1.4].
 However, the parameter surface of a minimal elliptic fibration is not uniquely
determined
by the birational equivalence class of the fibration:

\proclaim{\bf Proposition 0.4} Let $\bar \pi _i : \bar X_i \to \bar S_i, \
i=1,2$ be 2 birationally equivalent elliptic fibrations with $\bar X_i$ minimal
and $K_{\bar X_i} \equiv \bar \pi ^*(K_{\bar S_i} + \bar \bla _i),  \ \forall
i$. Denote by $ \varphi  : \bar S_1 \dasharrow \bar S_2 $ the corresponding
birational map.
  Then $(K_{\bar S_i} + \bar \bla _i) \cdot \Gamma =0$, for every curve
$\Gamma$ which is $\varphi$  or $(\varphi)^{-1}$-exceptional.
\endproclaim

\demo{\it Proof}
Let $q_i: Y \to \bar X_i$ the a resolution of the graph of the map $\mu :  \bar
X_1 \dasharrow \bar X_2$. Then $q_1 ^*(K_{\bar X_1}) = q_2 ^*(K_{\bar X_2})$
[Ko, 4.4].

Consider a curve $\Gamma$ $\varphi$-exceptional and any curve $C_1 \subset \bar
X_1$
dominating $\Gamma$. Let $\tilde C$ be the strict transform of $C_1$ on $Y$.
Then by commutativity of the diagram:
$$ 0 = (K_{S_2} + \bla _2) \cdot  ({\bar \pi _2 \cdot q_2})_* (\tilde C ) =
K_{\bar X _2}  \cdot ({q_2}_* \tilde C)  = K_{\bar X_1} \cdot C_1 = (K_{S_1} +
\bla) \cdot \Gamma .$$
\qed\enddemo

\medpagebreak

More generally the following holds:

\proclaim{\bf Propostion 0.5} Let  $\bar \pi _2 :\bar  X_2 \to \bar S_2$ be an
elliptic fibrations with $\bar X_2$ minimal.  Let $\bar \pi :  \bar X_1 \to
\bar S_1$ be a
minimal model  of $\bar X_2$ obtained as in (0.3).  Although the birational map
$\bar X_2 \dasharrow S_1$  is not in general an
 elliptic fibration, there exists a birational
 morphism $\tau :\bar S _1\to \bar S_3$ such that the induced birational
 map $\rho : \bar X_2 \to \bar S_3$ is not only a morphism but an elliptic
fibration.

Furthermore, for every curve $\Gamma$ which is $\tau$ exceptional, there is a
curve
 $C \subset \bar X_1$ dominating $\Gamma$ included on the exceptional locus of
the birational map $\bar X _1\dasharrow \bar X_2$ and additionally
 $ (K_{\bar S_1} + \bar \bla ) \cdot \Gamma =0$.
\endproclaim

\demo{\it Proof}  Then statement is obvious if $\bar S_2 \to \bar  S_1$ is a
morphism.  Otherwise consider the graph:
$$
\alignat3
& &&\ \ S &&\\
 & { f_1}  \swarrow && &&\searrow {f_2}\\
& \bar S_1 &&\overset{\varphi}\to{\dasharrow} &&\  \ \ \ \bar S_2
\endalignat
$$
Let $F \subset S$ be the union of the exceptional locus of $f_1$ and $f_2$.
Since the question is local we can assume $F$ connected.
Let $V \subset \bar S_3$ be an open set contaning $f_1(F)$.

Consider the set of curves $\{ C_i \} \subset \bar X_i$ which are included in
the exceptional locus of $\mu$ and dominate a curve $\Gamma \subset f_1 (F)$.
Then $ 0 =K_{\bar X_i} \cdot C_i = (K_{S_1} + \bla _1)  \cdot \Gamma$,  [Ko,
4.5].
Let $\tau : \bar S_1 \to \bar S_3$ the contraction of all such curve $\Gamma$
and
$$
\alignat3
& &&\ \ Y &&\\
 & {q_1}  \swarrow && &&\searrow {q_2}\\
& {\bar X_1} &&\overset{\mu}\to{\dasharrow} &&\  \ \ \ {\bar X_2}\\
& \downarrow \bar \pi _1 && && \bar \pi _2 \downarrow  \\
& {\bar S_1} &&\overset{\phi}\to{\dasharrow} &&\  \ \ \ {\bar S_2 }\\
& \tau \searrow && && \quad \\
& &&\ \ {\bar S_3} &&
\endalignat
$$
be a resolution of the graph of the map $\mu : {\bar X_1} \dasharrow {\bar
X_2}$.
Then
$$\rho (x) =
\left\{
\aligned  &\tau \cdot f_1 \cdot  f_2 ^{-1} (\bar \pi _2 (x)) \quad  \quad
\text{ if } \bar \pi _2 (x) \notin f_2 (F)  \\
 &\tau  \cdot \bar \pi _1  \cdot q_1 \cdot q_2 ^{-1}(x) \quad \text{ if } \bar
\pi _2 (x) \in   f_2 \cdot  f_1 ^{-1}
(V) \endaligned
\right.
$$
defines a rational elliptic fibration $\bar X_2 \dasharrow \bar S_3$.
We need  to show that $\rho$ is actually a  morphism.
Assume that $P \in \bar X_2$ is a fundamental point of $\rho$ and let
 $$
\alignat3
& &&\ \ W &&\\
 & { p_1}  \swarrow && &&\searrow {p_2}\\
& \bar X_2 &&\overset{\varphi}\to{\dasharrow} &&\  \ \ \ \bar S_3
\endalignat
$$
be the corresponding graph: note that
$p_1$ is a birational map.  The commutativity of the diargam and Zariski's Main
 Theorem imply that $p_1 ^{-1}(P)$ is connected and has dimension $0$, see for
example [Ha2, 5.2].  Then $\rho$ is actually a morphism.
\qed \enddemo

\bigpagebreak

 Since two-dimensional log terminal singularities are quotient [Ka4],  the
exceptional graphs
 on the minimal resolution  of $\bar S$ must be part of the Brieskorn's list of
surface quotient singularities [Br].

\bigpagebreak

If the fibration does not have multiple fibres singularites of $\bar S$ occur
when the $\Bbb Q$-divisor $\bda = \pi _* (K _{X/S})$ is not nef.

Because of the connection with the $C_{n,m}$ conjecture [Ka6] and the
birational classification of varieties, different notions of  ``positivity" of
the direct image of the relative dualizing sheaf have been introduced and
studied.
Here is a list of some ``positivity" results that are relevant to the case of
elliptic fibrations:
\proclaim{\bf Theorem   [Ka1] }
Let  $g: T \rightarrow W$ be any surjective morphism between non singular
projective varieties where the ramification divisor $\bsa$ has simple normal
crossings
and the local monodromies of  $ \Bbb R g_* \Bbb C _{X_0} $ are unipotent.
 Then $ g_* (K _{T/W}) $ is locally free and semipositive.
\endproclaim

\proclaim{\bf Theorem  [Kz]}
An algebraic fiber space with simple normal crossings branching has unipotent
monodromies if it is semistable in codimension 1.
\endproclaim

 These result show
 that all monodromies being unipotent,
or ``most" fibers reduced with simple normal crossings, is a sufficient
condition for $\pi _* K_{X/S}$ to be nef;
in the case of an elliptic fibration of a surface over a curve,   $ \pi _*
K_{X/S} = \bda$  has always non-negative degree.

Note that, in the case of invertible sheaves, semipositive is equivalent to nef
(see, for example [Ha1] and/or [Ka3]).

\bigpagebreak

\subheading{Notations}

Throughout this paper, $\pi :X \rightarrow S $  denotes an elliptic fibration
between smooth varieties, $\dim (X) =3$ and  $\bsa = \bsa _{X/S}$ a divisor
with simple normal crossings.

\noindent {\bf (0.4)} We assume  $S$ to be \it simple normal crossings minimal,
\rm that is, any contraction $\phi : S  \rightarrow S_1 $ of a log extremal
curve leads to $S_1$ singular or
$ \bsa _{X / {S_1}}$ not simple normal crossings.

\noindent {\bf (0.5)}
Let $E$ be an irreducible component of $\bsa$.
By abuse of notation we say that $E$ is of type $*$ if $\pi ^{-1} (P)$
is a fiber of Kodaira type $*$, for general $P \in E$.

\noindent {\bf (0.8)} A \it string \rm of lenght $r$ is a union $C=  \cup C_i$
of smooth rational curves $C_i$,
such that $C_i \cdot C_j = 1 $, for $|i-j| =1$ and $C_i \cdot C_j = 0 $, for
$|i-j| >1$.
A \it Hirzebruch-Jung  \rm string is a string in the above sense with $C_i ^2
\leq -2, \ \forall i$.

\noindent {\bf (0.9)} Consider a string of log extremal curves of length $r$.
Set  $S = S_0$  and denote by $\psi _i :  S_{i - 1} \rightarrow S_i$, $ 1 \leq
i \leq r$
the contraction of the log extremal curve $E_i \subset S _{i-1}$, and by
$\psi : S \rightarrow S _r $ the composition of the contraction morphisms.
Set $E  _{k,0}  = E _k$  and let  $ E   _{k,h} \subset S _{k-h-1}$  be the
strict transform of $E  _{k,h-1}$ .

The string is of type ${<} n , q {>}$, if
 ${<} n , q {>} $ represents the continued fraction:
 $$ \frac{n}{q} = b_1 - \cfrac 1\\
   b_2 -  \cfrac 1\\
   \dots - \cfrac 1\\
     b_r\endcfrac$$
and   $(E_1) ^2 =- b_1, \  (E_{2,1}) ^2 = - b_2 , \cdots  , (E  _{r,r-1})^2 = -
b_r$.

 For short we will write :
 $$ \frac{ f_k}{g_k} = b_k - \cfrac 1\\
   b_{k - 1} -  \cfrac 1\\
   \dots - \cfrac 1\\
     b_1\endcfrac$$  where $  f_k $ and $ g_k $ are
 positve integers defined by induction as follows:
 $$
\left\{
\aligned f_1 = b_1  \\
 g_1 = 1  \endaligned
\right.
 \qquad
\text{ and } \qquad
\left\{
\aligned b_k f _{k-1} - g _ {k-1} &= f_k\\
          f _ {k-1} &=  g _k.  \endaligned
\right.
\tag  0.10$$
In fact, $\frac{f _k}{g _k} = \frac{b_k f _{k-1} - g _ {k-1} }{ f _ {k-1} } =
b_k - \frac{g_{k-1}}{f_{k-1}} $.

Note that $\frac{g_r}{f_r} = \frac{q'}{n}$, where $ q' < n$ and $q q' \equiv 1
\ mod \ n$.
\medpagebreak

\noindent { \bf (0.11)} $K _k$ will denote $ K _{S _k}  $ , $\bda = \bda _0$,
  and
$ \bda _k$  will stand for $  ( \psi _{k-1})_* (\bda _{S _{k-1}})$.
Write
$$\bda _{k-1}  =\psi ^* _k  ( \bda _k ) + \alpha _k E _k ,  \quad \quad
K _{k - 1} = \psi ^* _k  (K _k) + \beta _k E _k ,$$ and for some rational
numbers $\alpha _k $
and $\beta _k$.

Thus:
 $$ \aligned  K _{k - 1} \cdot E _k
     &=  [ \psi ^* _k  ( K _ k  ) +  \beta_k  E _k] \cdot E_k
=    \beta_k  E_k ^2  \\
 \bda _{k-1}  \cdot E _k  &= [ \psi ^* _k  ( \bda _k ) + \alpha _k   E _k]
\cdot E_k
= \alpha _k   E_k ^2 \endaligned
$$

\medpagebreak

\noindent{ \bf (0.12)} If $S \rightarrow \bar S$ is a log contraction, with
$\bar S$ singular
and $S$ smooth, we say that $S$ is the \it log extremal resolution \rm of $\bar
S$.
Throughall this paper all the singular surfaces  $\bar S$ are obtained by
contracting log extremal curves and thus there exists a log extremal
resolution.

\bigpagebreak

\noindent{ \bf (0.13)} If the rational map $\bold J : S \rightarrow \Bbb P ^1$
induced by the elliptic fibration is constant along a divisor $E \in \bsa$,
with value $ \lambda$, we say that $E$ has $  J$-\it invariant \rm $\lambda$,
or, $J (E) = \lambda$.

\vskip 0.2in

\subheading{Properties}

The following is a list of first properties of log extremal curves (see the
proof of Theorem 0.3):

\proclaim{ Proposition 0.14}

\noindent {\bf (1)} If $ \bla \cdot E \geq 0 $, then $ (K _ S  + \bla ) \cdot E
< 0$
 if and only if  $ K _ S \cdot E < 0 $, that is $E$ contracts to a smooth
point.

\noindent {\bf (2)}
 Let $\bar E$ be a log extremal curve on a singular surface $\bar S$.
If $  \bar \bla \cdot \bar  E \geq 0 $, then the contraction of $\bar E$
``improves the singularity".

\noindent {\bf (3)}
$\bda \cdot E \in \Bbb Z$, for every curve $E$ on $S$, because $\bda$
is the divisor associated to the line bundle $\pi _* (K _{X/S})$.

\noindent {\bf (4)} If  $\psi : S \rightarrow S_1 $ is a log-contraction to a
smooth point,
then $( K _{ S _1}  +   \bda _1 ) $ is  a Cartier divisor;
consequently $\bda _1 \cdot E $ is an integer
 for any curve $ E $ on $\bar S$.
\endproclaim

\demo{ \it Proof}
(1),(3), (4) follow from 0.3.

(2) Let $ \phi: S \rightarrow \bar S$ be the minimal log extremal resolution
(see 0.12).
If $ E$ denotes the strict transform of $\bar E$, then the assumptions imply  $
K _ S \cdot  E  < K _ {\bar S} \cdot  {\bar E}  <0$.

Thus $\phi$ can  can be re-factored, first contracting  $ E$ to a smooth point.
\enddemo

\vskip 0.4in

\heading  {\bf 1.  Strategy }
\endheading

 In \S 2 and \S 4 we  investigate the singularities arising from contractions
of
 log extremal curves in $\bsa - \bsa ^m$.

\medpagebreak

{}From now on,  $\phi : S \rightarrow \bar S$ denotes the contraction
of all the log exceptional curves $\Gamma$ in $\bsa - \bsa ^m$; then
 $ (K_S + \bla ) \cdot \Gamma = (K_S + \bda ) \cdot \Gamma$.
 The key point in the analysis of such contractions is that $\bla \cdot \Gamma
\geq \bda \cdot \Gamma$ and that $\bda \cdot \Gamma$ has to be an integer
[0.14, part 3].

Our strategy is to factor the contraction morphism $\phi: S  \rightarrow \bar S
$ into 3 morphisms $\phi = \varphi _1 \cdot \varphi _2 \cdot \varphi _3 $:
$$ \varphi _1  : S \rightarrow S_1 \text{  ,  } \varphi _2 : S _1 \rightarrow T
\text{   ,  } \varphi _3 : T \rightarrow \bar S \text{  .}$$

\roster
\item  $ \varphi _1$ denotes the composition of contractions of successively
 log extremal curves to smooth points and $S_1$ is a smooth surface. Set $\bda
_ 1 = {\varphi _1} _* (\bda)$.
\item  $\varphi _2$ is the morphism associated to the contraction of log
extremal Hirzebruch-Jung  curves. $T $ is a singular surface and
 $S_1$ is its minimal resolution. Set $\bda _ T = {\varphi  _2} _* (\bda _1)$.

\item $ \varphi _3$ denotes the contraction of all the remaining log extremal
curves. Then  $ \bar \bda =( \varphi _3) * (\bda _T) = \phi _* (\bda )$.
 \endroster

The order of contraction of log extremal curves is arbitrary: this one has been
chosen because of its logical convenience.
The argument is nevertheless general,
because the minimal surface (with log terminal singularities) associated to
$S$ and $\bla$ [Sa, Thm 1.4].

First we will set up a combinatoric contraction algorithm
for (1) and (2)  (in 2.2., 2.5.) and then prove that all the curves contracted
by $ \varphi _3$
have self-intersection  $(-1)$ on the minimal resolution $S_1$.
We will show that if $E_i$ is contracted by $\varphi _3 $ , then   $( K_ {T} +
\bda  _{T}) \cdot \hat E _i = 0$, where $\hat E_i$ is the strict transform of
$E_i$ on $S_1$.
These two properties will allow us to re-factor $\varphi _2
\cdot \varphi _3$, first contracting $\hat E_i$ to a smooth surface $S' _{i_1}$
and then applying the previous algorithm to the remaining contractions.
This process thus produces an iterative contraction algorithm; eventually
$\phi$ is factored as:
$$ \varphi'_1 : S \rightarrow S' \text{  and  } \varphi'  _2 : S' \rightarrow
S'_2 =
\bar S \text{  ,}$$
where $S'$ is the smooth minimal resolution of $ \bar S$ , $ \bsa _{X/ S'}$ is
not a divisor with simple normal crossings (if $S \neq S'$) and $\varphi' _2$
 contracts only Hirzebruch-Jung  strings.

Note that $\varphi ' _1$ is not a log extremal contraction.

\smallpagebreak

The contraction of $\varphi_1$ (respectively $\varphi_2  , \ \varphi_3$) will
be denoted by the first (respectively second and third) step.

\vskip 0.2in

In \S 2 and \S 4 we state  and prove the major points in the algorithm;  the
technical details are treated in
\S 5.

\vskip 0.4in

\heading{\bf 2. The first two steps}
\endheading

First we consider the case of a log extremal curve on a smooth surface
 with $J$ invariant $ \infty$.

The following proposition says that log extremal strings  with $J = \infty$ do
not contract to singular points of $\bar S$.

Proposition 2.2., 2.1 and 2.5. together form the basic for the combinatorics of
our contraction algorithm.
 In fact, 2.5. states that a log extremal curve $E_1$ on a smooth surface must
have intersection $-1$ with $ (K _S + \bda)$ while  2.5. gives bounds on the
self-intersection number
of $E_1$ and on the sum of the coefficients (in $\bda$)
of the intersecting curves.

The second parts of 2.2. and 2.5.  tell us how to continue the contraction
algorithm in (1) and (2).

\medpagebreak

\proclaim { \bf Proposition 2.1} Let $\pi : X \rightarrow S$ be an elliptic
fibration such that the $ J$ invariant function extend to a morphism $\bold J:
S \rightarrow \Bbb P ^1$. Assume that  $\bold J$ has a pole over a log extremal
curve $B _1 $ and denote by $\psi _1 : S \rightarrow S_1$ the corresponding
contraction.

\noindent {\bf (a)} $\bold J:  S_1 \rightarrow \Bbb P ^1$ is a morphism
$\psi _1^* (J^1 _\infty) = J_\infty$, where $ J^1 _\infty$ denotes the divisor
of poles of $\bold J$ on $S_1$.

\noindent {\bf (b)} If $S$ is smooth and $\bda \cdot C \in \Bbb Z$, for all
1-cycle $C$ on $S$, then $(B_1)^2 \geq -2$. $(B_1)^2 = -2$ only if the fiber
over the general point of $B_1$ is of type $I_a ^*$.
If $(B_1)^2 = -1 $, then $\bda_1 \cdot C \in \Bbb Z$, for all 1-cycle $C$ on
$S_1$.

\noindent {\bf (c)} Let $B_1$ and $S_1$ be as in (b), with $(B_1)^2 = -2$. Let
$B_{2,1}  \in \bsa$ be another curve such that $B_{2,1} \cdot B_1 = 1$  and
$B_2 = (\psi _1) _* (B_{2,1})$ is a log extremal curve on $S_1$.  Denote by
$\psi _2 : S_1 \rightarrow S_2$ the corresponding contraction.  Then $S_2$ is
smooth  and $(K_S + \bda) \cdot B_{2,1} =0$.

\endproclaim

Recall that is $S$ is the smooth surface image of log extremal contractions,
then   $\bda \cdot C \in \Bbb Z$, for all 1-cycle $C$ on $S$. Note  also that
this agrees with [Kz].

\demo{\it Proof}

\noindent {\bf (a)} Consider $H$ a smooth irreducible element of $  J _
\infty$.
Since $H$ is a fiber of $\bold J$, $H \cdot B _1 = 0 $. Thus the statement
follows directly from the definition.

\noindent {\bf (b)} Write
$$\bda = \sum a_i D_i + \delta (B_1) B_1 + (1/12)J_\infty, $$
where $D_i \neq B_1 , \forall i$; $\delta (B_1) = 1/2$ when
$B_1$ is of type $I_a ^*$ or $\delta (B_1) = 0$ when $B_1$ is of type $I_a$
[Kd].

 Set $\bda _ 1 = \psi _* (\bda)$ and $\bar D _i = \psi _* (D_i)$.  Thus $\bda
_1 = \sum a_i \bar D_i + (1/12) J _\infty$. We will prove the statement
comparing $\bda$ and $\psi _1 ^* (\bda _1)$.  Set:
$$\aligned
\bda &= \psi _1 ^*(\bda _1) + \alpha_1 B_1, \ \alpha_1  \geq 0  \\
D_i &= \psi _1 ^*(\bar D_i) - \gamma_i  B_1, \ \gamma_i  \geq 0
\endaligned
$$
It follows that $ (1/2) \delta = \sum a_i \gamma _i + \alpha _1$ and $ 0  \leq
\alpha  _1\leq 1/2$, with $\alpha _1 = 1/2$ only when $\delta = 1/2$ and
$\gamma _i = 0 , \forall i$.
Then (5.5) applies and $1/2 \geq \alpha = 1 + \frac{1}{(B_1) ^2}$. This shows
that
$ (B_1) ^2 \geq  -2$ and $ (B_1) ^2 =  -2$ only when $B_1$ is of type $I_a ^*$
and $\gamma _1 = 0 , \forall i$.

\noindent {\bf (c)} Set $B_{2,1} ^2 = -b_2$.
Then
$$ \bda_1 \cdot B_2 = \bda \cdot  B_{2,1} -1/2 \geq 1- 1/2 = 1/2.$$
Since $K_S  = \psi ^*(K_{S_1})$ and $B_2$ is log extremal, the above inequality
implies
$$0> (K_{S_1} + \bda _1) \cdot B_2 = b_2  -2 + \bda _1 \cdot B_2 \geq b_2
-3/2.$$
Thus $b_2=1$; furthermore $\bda_1 \cdot B_{2,1} =1$ and  $(K_S + \bda) \cdot
B_{2,1} =0$.\qed
\enddemo

\bigpagebreak

\proclaim{ Proposition  2.2}
Let $\{ E_1 , E_{2,1} , \dots , E_{r, r-1} \}$  be a string of log extremal
curves on a smooth surface $S $.
($\bsa$ is not necessarily simple normal crossing.)
 Assume that  $\bda \cdot E_{j, j-1}  $  is an integer,
  $\forall j $ and that the string is represented
by ${<} n , q {>}$.
Then in the notation of  (0.9) we have:

\noindent {\bf (1)} $(K_ S + \bda )\cdot E_1 = -1$
 and $(K_S + \bda) \cdot E _{j,j-1} = 0$ , $\forall j, \  2 \leq j \leq r$.

\noindent {\bf (2)}  $( K _S  +  \bda ) \cdot \bda = ( K _{S _r} +  \bda _r  )
      \cdot \bda _r  + \frac{1}{n} -1 $.

\endproclaim

\demo{ \it Proof of Proposition 2.2.}

\noindent {\bf (1)} By (5.1), (5.2), (5.5) and (0.4) we have
$$( K _S + \bda ) \cdot E _ 1 =
(\alpha _1 + \beta _ 1) (E _1) ^2, $$
and thus  $\frac{1}{f_1} \ \frac{- f_1}{g_1}= - 1 .$

Similarly, for any $ k \geq 2$ we can write:
$$( K _S + \bda ) \cdot E  _{k , k -1} = (\alpha _k + \beta _ k) (E _k) ^2 +
\alpha _ {k-1} + \beta _{k-1} = \frac{ - f _k }{f_k g_k} +
\frac{1}{  f _{k-1}}= 0 .$$

\smallpagebreak

\noindent {\bf (2)}  By  (0.11) , $( K _S + \bda ) \cdot \bda =
 ( K _{ S_r} + \bda _r  ) \cdot  \bda  _r  +
\sum \Sb k\endSb
\alpha _k (\alpha _k + \beta _ k) (E _k) ^2  \text{  .}$

\noindent Substituting the expressions for $\alpha_k $ and  $\beta_k$   derived
in
5.3. and 5.5. we have
$$\sum _{k=1} ^r
\alpha _k (\alpha _k + \beta _ k) (E _k) ^2 = \sum _{k=1} ^r
  \frac{ f _ k - g _ k }{  f _ k} \cdot \frac{1}{f_k} \cdot
\frac{ - f _ k }{ g _ k } =  \sum _{k=1} ^r
- 1 / g _ k + 1 / f _ k.$$
By  (0.11) $\sum \Sb k\endSb
(- 1 / g _ k + 1 / f _ k )=
- ( 1- 1/n)$ and thus
$$ ( K _S + \bda ) \cdot \bda =
( K _{ S_r} +  \bda _r ) \cdot  \bda _r - ( 1- 1/n) .
 \qed $$
\enddemo

The following corollary is a particular case of Theorem 2.5.: however, it is
also an immediate application of 2.2.

\proclaim{\bf Corollary 2.3}
Let $ \{ E_1, \dots , E_ {r,r-1} \} $ be a Hirzebruch-Jung string of log
extremal
curves on a smooth surface $S$ represented by ${<} n , q {>}.$
 Assume that the string is a connected component of $\bsa$.
 Then:

\noindent {\bf (1)}
$\bda = (1 - \frac{g_r}{f_r} ) E_{r, r-1} + \dots (1 - \frac{g_2}{f_r} ) E_{2,
2-1}  + (1 - \frac{g_1}{f_r} ) E_1. $

\noindent {\bf (2)}
 $ r \leq 5$; $ {<} n , q {>}$ and $ ( a_1, \dots, a_ r )$  (the coefficients
of $\{ E_1, \dots E_ {r,r-1} \}$ in $\bda $)  are described in the table below.
$$
\vmatrix
 n & & q & & a_1  & & a _2 & & a_3 & & a _ 4 & & a _5 \\
2 & & 1 & & 1 / 2 \\
3 & & 1 & & 2 /3  \\
6 & & 1 & & 5 /6 \\
4 & & 1 & & 3/4 \\
4 & & 3 & & 3/4  & & 1 / 2 & & 1 /4 \\
3 & & 2 & & 2 /3  & & 1 / 3 \\
6 & & 5 & & 5 /6 & & 2 /3  & & 1 / 2 & & 1 / 3 & & 1/6 \\
\endvmatrix
$$
\endproclaim

\demo{\it Proof}
Part (1) follows from 2.2 and 5.5.
It follows also that  $a _r = \alpha _r = 1 - \frac{g_r}{f_r} = 1 - \frac{
q'}{n}$, where $ q' < n$ and $q q' \equiv 1 \ mod \ n$.
Note that $ \{ E_1, \dots , E_ {r,r-1} \} \subset \bsa $ and  that $p$ singular
implies $ \bold J (E_j ) \neq \infty $, for every $j$ (see 2.1 ).
Hence $ a _r \in \{1 / 2,\ 2 /3 , \ 5 /6,\  3/4, \ 1 /4, \ 1 / 3 ,\ 1 /6 \}$
[Kd]. \qed
\enddemo

\vskip 0.2in
2.3 describes the log extremal Hirzebruch-Jung  strings which are connected
components
of $\bsa _{X/S}$.  The last part of this section contains numerical conditions
which, together with  2.1 and 2.2., describe the contraction of any log
extremal string.

\proclaim{ \bf Proposition 2.4}
\roster
\item Let $E_1$ be a log extremal curve  of self-intersections
$- b _ 1 $.

\noindent Write $\bda =  e_1 E_1 + \sum c_i C_i$, where each $C_i$ is an
irreducible curve,  $c_i$ and $e_1$
are the rational number described in $0.2$.

 If $(K _S + \bda) \cdot E _1 = -1$, then
$$0 \leq \sum\Sb i \endSb  c_i (C_i \cdot E _1) = 1 - ( 1- e_1) b_1 < e_1$$
In particular, $ b _1 \leq ( 1- e_1) ^{-1} $.
\item Let $ E _{2,1}$ be another log extremal curve such that $ E_1 \cdot E
_{2,1} = 1$ and $(E_{2,1} ) ^2 = -b_2$ (the notation is as in 0.9).

\noindent Write $\bda =  e_1 E_1 + e_2 E_2 +  \sum c_i C_i$.
If $(K _S + \bda) \cdot E _{2,1} = 0$ , then
$$ 0 \leq \sum\Sb i \endSb  c_i (C_i \cdot E _{2,1}) = 2 - ( 1- e_2) b_2 -
e_1$$
In particular, $ b _2 \leq \frac{ 2- e_1}{  1- e_2} $.
\endroster
\endproclaim

\demo{\it Proof} $E_1$ and $E_{2,1}$ are rational curves, hence:
 $ \bda \cdot E_1 = ( K_S + \bda ) \cdot E_1 - b_1 +2 $;
 $ \bda \cdot E_{2,1} = ( K_S + \bda ) \cdot E_{2,1} - b_2 +2 $.
\roster
\item $(K_S + \bda) \cdot E _1 = -1
\Longleftrightarrow
 \bda \cdot E_ 1= 1 - b _ 1 $
thus $$\aligned
 1 - b _ 1 & =  - e_1 b _ 1 + \sum\Sb i \endSb c _ i (C _ i \cdot E _1) \quad
\text{ i.e. } \\
 1 - ( 1- e_1) b _ 1 &=\sum \Sb i \endSb c _ i (C _ i \cdot E _1)
\geq 0.\endaligned
$$

In particular  $b _ 1 \leq ( 1 - e_1) ^ {-1}$.
\item Similarly,
$$ \aligned e _1 - e_2 b_2 + \sum\Sb i \endSb c_i (C _i \cdot E _{2,1})
&= \bda \cdot E _2 = 2 - b _ 2  \quad \text{ i.e }\\
 2 - (1 - e_2) b _ 2  - e_1 &=\sum \Sb i \endSb c _ i  (C _i \cdot E
_{2,1})\geq 0. \endaligned
$$
In particular $ b _2 \leq \frac{ 2- e_1}{ 1- e_2} $.\qed
\endroster
\qed
\enddemo

\proclaim { \bf Proposition 2.5}

Let $\{ E_1, \ E_{2,1} \}$  be two
log extremal curves on the smooth surface $S$.  As in (0.9) we denote by $\psi
_1 : S \rightarrow S_1$ the contraction of $E_1$ and by $\psi _2 : S _1
\rightarrow S_2$ the contraction of $E_2$ .
If $S_1$ is smooth so is $S_2$.
\endproclaim

\demo{\it Proof}
By 2.1  we can assume $J (E_1) \neq \infty .$

Denote by $e_1$ (\it resp \rm $e_2$) the coefficient of $E_1$
 (\it resp \rm $E_{1,2}$)  in $\bda$ and use the same notation of 2.4.
Since  $E_1$ contracts to a smooth point $- (b_2 -1) = (E_2)^2$.

 Assume that $b_2  - 1 \geq 2$ and derive a contradiction.
 By 2.2 $(K_{S_1} +\bda _1 )\cdot E_2 =-1 $ and
$$ \bda _1 \cdot E_2 = -1 -K_{S_1} \cdot E_2 = -1 +2 - (b_2 -1) = 1- (b_2 -1)
.$$
On the other hand, since $E_1$ contracts to a smooth point
 $\bda \cdot E_1 =0$ and
$ \bda _1 \cdot E_2 \geq e_1 - e_2 - (b_2 -1) e_2$.
Note also that $e_1 - e_2 \geq 0$. Combining the two expressions for $\bda
_1 \cdot E_1$ we have:
$$1- (b_2 -1)  \geq e_1 - e_2 - (b_2 -1) e_2, \text{ that is}$$
$e_1 - e_2 \leq 1-(b_2 -1) ( 1- e_2)$.
 If $b_2 -1 \geq 2$, then $e_1 - e_2 \leq 2 e_2 -1$ and thus $ e_1 \leq 3e_2
-1$.
Recall that $S$ is assumed to be simple normal crossing minimal (0.4), and thus
$e_1 - e_3 \geq 1/3$.  Combining the two inequalities we get $e_1 \leq 3 e_1
-2$
or $ 1 \leq e_1$, which is a contradiction.
 \qed
\enddemo
Proposition 2.5, which depends on 2.1,  also implies 2.1 as a particular case.

\vskip 0.3in

Now we are ready to  apply the previous computation in order to describe
the log extremal curves contracted by $\varphi _1$ and $\varphi _2$.

Recall that:
\noindent

\noindent{\bf (1) First Step}
 $\varphi _1 : S \rightarrow S_1$ is the contraction of all log extremal curves
to smooth points.

\noindent{\bf (2) Second Step}
$\varphi _2 : S_1 \rightarrow T$  contracts all the chains of disjoint
Hirzebruch-Jung  curves.

\noindent (The case when $J $ has poles on the  whole string  is considered in
(2.1); thus we can assume that the first curve to contract has $J$ invariant
$\neq \infty $; furthermore, the first curve to contract is disjoint from the
image of the previous exceptional smooth strings by Prop. 2.5.)
$S_1$ is the minimal resolution of $T $.

\medpagebreak

\proclaim{ Theorem 2.6}

\noindent After the first step (1) the possible singularities of $\bsa _{X/T} $
coming from contractions of curves with $J \neq \infty$ are the one listed in
Table 1.

\noindent After the second step (2) the Hirzebruch-Jung strings
which can occur are of type ${<} n , q {>}$ , $ n \leq 6$.:

Furthermore the type of the corresponding singular fiber are
 listed in Table 2, 3, 5.
\endproclaim

\demo{\it Proof} If $J(E_1) = \infty$ then the statement follows form (2.1).
Assume  $J(E_1) \neq \infty$.
The application of the contraction algorithm is straightforward and  here we
work out only the particular case of $e_1 = 2/3$. Further details are avaliable
upon request.

Recall that we assume $S$ to be simple normal crossing minimal.
The notation is as in 2.4

\smallpagebreak

\noindent {\bf Case 1:} If $ b_1 = 1$, then $(K_S + \bda) \cdot E_1 = -1$
and $\sum c_i = 2/3$ (2.4).  Thus $\{ c_1, c_2, c_3 \} = \{ \frac{1}{3},
\frac{1}{6}, \frac{1}{6}\}$.
Again by 2.4 only $C_1$ can be log extremal.  Set $C_1= E_{2,1}$.  Then $0 =
E_{2,1} \cdot \bda _1 = \frac{1}{3} + \frac{1}{6} + \frac{1}{6}$.

\smallpagebreak

\noindent {\bf Case 2:} $ 2 \leq b_1$. By 2.4.,  $b_1 \leq 3$.

{\bf Case $2_A$} $b_1 = 3$. The curve $E_1$ (type $IV ^*$) has to be a
connected component of $\bsa _{X/ {S_1}}$.

{\bf Case $ 2_B$ } $b_1 = 2$, $\sum c_i = 1/3$.
Then $ c_1  = 1/3$ (type $IV$), or $\{ c_1, c_2 \}=  \{  \frac{1}{6},
\frac{1}{6}\}$ (types $II , II$).
$C_1$ can be part of a log extremal Hirzebruch-Jung string only if it is of
type $IV$ (2.4, part 2
and 2.2).  Set $C_1 = E_ {2,1}$.  Then $b_2 = 2$ and  $\{ E_1 , E_{2,1} \}$ is
a connected component of  $\bsa _{X / {S_1}}$.
\qed
\enddemo

\medpagebreak

 In the following pages $r$ denotes the length of the string.
The roman numbers denote the type of the singular fiber in the Kodaira
 classification, while the arabic ones, between parenthesis,  are the
self-intersection numbers of the exceptional curves.

 The solid (resp. empty) circles represent the curves in
$\bsa _{X/T}$ which are (resp. are not) log extremal.

Note that except the $<5,2>$ and $ <5,3> $ cases (where the exceptional curves
have selfintersection (-2,-3 ) and (-3,-2)), all the strings
of length $ \geq 2$ are of $(-2)$ curves.

\newpage

\pageno=17

\heading{\bf 3.  Examples}
\endheading

 \medpagebreak

\example{Example 3.1 - 3.7}

In the 2 following examples,
$\zeta = e ^ {2 \pi i / n}$ is a root of unity and
 $E =  { \Bbb C } /{{<} 1, \zeta {>}}$ a smooth elliptic curve.  Consider $\bar
X=
(\Bbb C ^2 \times E) / G $, where $G = {\Bbb Z} / {n \Bbb Z}$ and the
action is specified below.
In (1) and (2) $\bar S= \Bbb C ^2  / G$ has a singular point at the origin and
$\bar \pi : \bar X \rightarrow  \bar S$ is an elliptic threefold with terminal
singularities.
Let $\phi  : S \rightarrow \bar S$ be the minimal resolution of $\bar S$
 and $\pi : X \rightarrow S$ the resolution of the induced fibration:
$$
\CD
X         @.         \bar X \\
@V{\pi }VV              @V{\bar \pi}VV  \\
S        @>\phi>>      \bar S \\
\endCD
$$
\roster
\item Let  $n=2,3,4,6$  and $G = {\Bbb Z} / {n \Bbb Z}$
 act on  $(z_1, z_2, e) \in ( \Bbb C ^2 \times E)$ as follows:
 $ (z_1, z_2, e) \mapsto (\zeta z_1,\zeta  z_2,\zeta ^{-1} e)$.

 The exceptional graphs of $\phi$ are of type $ {<} n, 1 {>}$ and  $\pi ^{-1}
(P)$ is a singular
fiber with monodromy eingenvalue $e^{{2 \pi i (n-i)} / n}$ (corresponding
respectively to the Kodaira types $I_0^*$, $IV^*$, $III^*$, $II^*$).

\item
Let $n=3,4,6$ and $G = {\Bbb Z} / {n \Bbb Z}$ act on $(z_1, z_2, e) \in
 ( \Bbb C ^2 \times E)$ as: $ (z_1, z_2, e) \mapsto (\zeta z_1,\zeta  ^{-1}
z_2, \zeta e).$

 The exceptional graphs of $\phi$ are of type  $ {<} n, n-1 {>}$ and the
corresponding divisors of singular fibers are:  $IV^*$-$IV \ (n=3)$;
 $III^*$-$ I_0 ^*$-$ III \ (n=4)$ and $II^*$-$ IV ^*$-$I_0 ^*$-$IV$-$ II  \
(n=6)$.
\endroster

\endexample

\smallpagebreak

\example{Example 3.8: $q=1$, $n=5$}

Let $\Gamma$ be the rational normal curve in $\Bbb P ^5$ given by the embedding
$[x,y] \mapsto [ y^5,  y^4 x ,   y^3 x ^2 , y^2 x  ^3 , y x ^4,  x^5]$.
Let $T$ be the  projective cone over $\Gamma$ and $\tilde T$ its minimal
resolution with exceptional divisor $D$.

The automorphism $g:$
$[x,y] \mapsto  [x, \zeta y] $  of $\Bbb P ^1$ (where  $\zeta = e ^ {2 \pi i /
6}$) induces the automorphism of the cone:
$$   ( z_1,  z_2,   z_3  , z_4 , z _5,  z_6) \mapsto (\zeta ^5 z_1, \zeta ^4
z_2  ,  \zeta ^3 z_3  , \zeta ^2 z_4 , \zeta z_5 , z_6 ) .$$
It easy to see that $L = \{ z_j =0 \}, j=1,2,3,4,5$  is a line entirely fixed
by $G$, the group generated by $g$;  $D$ is invariant under $G$ and the  unique
fixed point $P \in D$  generates a singularity on  $ \tilde  T  / { G} $.

\midinsert
\vspace{3in}
\endinsert

Let $\rho : \tilde S \rightarrow \tilde  T$ be the blow up at $P$,  with
exceptional divisor $\tilde F$ ($\tilde F ^ {2} = -6$); $\tilde D$ the strict
trasform of $D$ has self-intersection $-6$.  $\tilde S /{G}$ is  then a smooth
surface and ${\tilde D  }/ G $ can be contracted to a smooth point $\chi :
\tilde S  \rightarrow S$; furthermore $ F^2 = {\chi _*  (\tilde F) }^2 = -5$.

\medpagebreak

We can extend the action of $g$ to $T \times E$
  by sending $e \mapsto \zeta e$.  The ramification divisor of
the induced fibration $\pi : X \rightarrow S$  is supported on $F$ and $L$.
  A computation in local coordinates shows that the singular fiber are
 respectively of type $II^*$ and $II$ and that $X$ has terminal singularities.

\endexample

\vskip 0.3in

\heading{\bf 4.  More contractions}
\endheading

 Now, the particular goal is to investigate the case of new contractions, and
in particular look for new exceptional graphs.
{}From now on, $E$ is a  log extremal  curve on $T$ and  $\hat E$ the strict
transform of $E$ on $S_1$.

 The curves that can contract now (and could not have been contracted with the
first step $\varphi _1 : S  \rightarrow S_1$ and
the second step $ \varphi _2  : S_1  \rightarrow T$) belong to the following
classes:

(i) curves $ E \notin \bsa _{X/T}$  whose image pass
through the singularities of $ T $,

(ii) curves intersecting 1 or more
exceptional graphs, but not part of a Hirzebruch-Jung string.
[These curves are represented by empty cirles in the previous pages.]

\smallpagebreak

The following is an example of case (i).

\medpagebreak

\example{Example}
Let $\pi _2 : X_2 \rightarrow S _2$ be an elliptic fibration between smooth
varieties ($\dim (X) =3$); let $C \in \bsa _{X_2 / S_2}$ such that $\pi _2
^{-1} (p)$ is a singular fiber of type $I_0 ^*$, for general $p \in C$. [In a
neighborood of $p$ the local equation for $X_2$ is given by  the
desingularization of the Weierstrass equation: $y^2 = x^3 + s^2 t^2 x + s^2 t^2
 $.]

 Let $\psi _2 : S_1 \rightarrow S_2$ be the blowup of $S_2$ at $p$ with
exceptional divisor $E_2$. Denote by $\pi _1 : X_1 \rightarrow S_1$ the
resolution of the induced fibration. The singular fiber over the general point
of $E_2$ is again of type $I_0 ^*$.

 Let $\psi _1 : S_0 \rightarrow S_1$ be the blowup of $S_1$ at the intersection
point of $E_2$ with the strict trasform of $C$; let $E_1$  be the exceptional
divisor and  $E_{2,1}$ the strict trasform of $E_2$.  Denote by $\pi _0 : X_0
\rightarrow S_0$ the resolution of the induced fibration. The fiber over the
general point of $E_1$ is a smooth elliptic curve.

If $K_ {X_0 }= \pi ^* (K_{S_0} + \bda _0)$ (this happens if $X_0$ is the
resolution of a Weierstrass model), then $K_ {X_j }= \pi ^* (K_{S_j} + \bda _j)
\ \forall j= 1,2$,
$\bla _1 = \psi  _2 ^* (\bla _2 )$, $\bla _0 = \psi  _1 ^* (\bla _1 ) -E _1$
and thus
$ (K_{S_0}+ \bla _0 )  \cdot E_{2,1} = -1 $ while $ (K_{S_0}+ \bla _0 )  \cdot
E_1 = 0 $.

Thus there exists a  first log extremal contraction $\varphi _1$ contracting
$E_{2,1}$ to a singular point and then a second log extremal contraction
$\varphi _2$, which contracts $E$, the image of $E_1$, to the smooth point $p
\in S$.
$\varphi _1 \cdot \varphi _2$ is a log extremal  re-factorization of $\psi  _1
\cdot \psi  _2$.
\qed
\endexample

We will show that all the new contraction are of this type; in particular, we
 prove that if $E_i$ is contracted by $\varphi _3 $ , then   $( K_ {T} + \bda
_{T}) \cdot \hat E _i = 0$, where $\hat E_i$ is the strict transform of $E_i$.
In particular we will show:

\proclaim{ Proposition 4.1}
If  a curve $E $ is log extremal on $T$, then $E$ passes through at most a
singular point of $T$.
In particular the contraction of $E$ improves the singularities and $({ K_{S_1}
+ \bda _1 }) \cdot \hat E = 0$.
\endproclaim

This  property will allow us to re-factor $\varphi _2
\cdot \varphi _3$, first contracting $\hat E_i$ to a smooth surface $S' _{i_1}$
so that the previous algorithm (2.2 and 2.4) applies.
In fact:
\proclaim\nofrills{\quad \usualspace}
 If $E$  is  a log extremal curve on  a smooth surface $ V$ such
that
  $(E) ^2 = -1$ and $( K _V  +  \bda  ) \cdot E = 0$, then
 $$( K _V  +  \bda  )  = \psi ^* (K _{ S_ 1 }  +  \bda _1 ) \text{ ,} $$ where
$\psi  : S \rightarrow   S _1$ is the contraction of $E$.

In particular if $( K _S  +  \bda  )  \cdot C \in \Bbb Z$,
for a  1-cycle $C$ on $S$, also
$(K _{ S_ 1 }  +  \bda _1 ) \cdot \psi _* (C) \in \Bbb Z$.
\endproclaim

This process thus produces an iterative contraction algorithm; eventually
$\phi$ is factored as:
$$ \varphi' : S \rightarrow S' \text{  and  } \varphi'  _2 : S' \rightarrow
S'_2 =
\bar S \text{  ,}$$
where $S'$ is the smooth minimal resolution of $ \bar S$ , $ \bsa _{X/ S'}$ is
not a divisor with simple normal crossings (if $S \neq S'$) and $\varphi' _2$
contracts only Hirzebruch-Jung  strings.

Again  most of the technical details are relegated in \S 5.

We start with the following:

\definition{ Definition 4.2}
Let $  \Cal E _1 = \{ E _{1_1}, \dots E_{r_1, r_1 -1} \}   \dots  \Cal E_s = \{
E _{1_s}, \dots E_{r_s  ,r_s -1} \}$
be $s$ disjoint Hirzebruch-Jung  strings, as in (0.9), and F a smooth rational
curve
such that

\noindent $ E_{r_j -1,r_j} \cdot F = 1$ , $\forall j $, $ 1 \leq j \leq s $ .

Then $F$ is said to be at the  {\rm \bf feet} of the strings
$ \Cal E_1  \cdots \Cal E_s .$

\enddefinition

By inspection on the graphs in Table 2, 3 and 4,
and from general consideration, we deduce  the following Remarks:

\proclaim{ Remarks 4.3}
 \roster
\item If $E \notin \bsa _{X/T}$, then $(\hat E) ^2 = -1$ (0.14, part 2).

\item If  $ E \in \bsa _{X/T}$ , then it is always at the ``feet" except in the
cases
at the botton line of Table 3 ($ II ^*, \ I_0 ^*$), and a case in Table 4 ($II
^*, \ IV^*, \ IV$).

\item If $ E$  intersects only one
exceptional graph  at the ``foot", then, by construction, $(\hat E) ^2 = -1$.

\item If $( \hat E) ^2 = -1 $, and $ \hat E$ intersects
  some chain of (-2) curves $\Cal  E$, then  $\{ \hat E, \Cal E \}$ contracts
to a smooth point.

\item  If $\hat E$ intersects 2 or more chains of (-2) curves, then  $( \hat E)
^2 \leq (-3)$.

\item  $E$ cannot be at the ``feet" of 4 or more exceptional curves
      and being exceptional on the minimal resolution ($\bar S$ has at most
quotient singularities). \rm

\endroster
\endproclaim

 In particular, we need to investigate the case of log extremal curves
intersecting
 different chains.
A crucial ingredient in the proofs is the following:

\proclaim{ Lemma  5.6}
Let $\{ \Cal E_1,\cdots  ,  \Cal E_s \}$
be $s =2, 3$ disjoint Hirzebruch-Jung  strings, and $F_0$ a smooth rational
curve
at the feet of $ \{ \Cal E _1,\cdots  ,\Cal E_s \}$.

Then the image of $F_0$ is not log extremal on $T$.
\endproclaim

Now we consider a log extremal curve $E \notin \bsa _{X/T}$.

\proclaim{\bf  Remark 4.4}
Lemma 5.6 implies that if $\hat E$ intersect more than 1 Hirzebruch string,
then $\hat E$ must also intersect a Hirzebruch-Jung string  $\Cal E$of lenght
$r \geq 2$;
only the following cases can occur:

\noindent {\bf (A)}  $\hat E$ intersects at most 1 Hirzebruch-Jung string

\noindent {\bf (B)}  $\Cal E$ is of type $<5,2>$ and $\hat E$ intersects 1
other curve of self-intersection $\leq -3$

\noindent {\bf (C)} $\Cal E$ is of type $<5,2>$ and $\hat E$ intersects 1 other
curve $C$ of self-intersection $ -2$

\noindent {\bf (D)} $\Cal E$ is of type $<5,3>$ and $\hat E$ intersects 1 other
curve $C$ of self-intersection $ \leq -3$

\endproclaim

\bigpagebreak

\demo{\it Proof of 4.1}

Case (1): $E \in \bsa _{X/S}$. From Remark 4.4. and  Lemma 5.5. follows that
$\hat E$ has to be at the feet of only 1 string of log extremal curves or $\hat
E$ has to be the curve of type $II$
(as in 4.3, part 2).

 In the first case the statement follows from  4.3 (part 3) and 5.10.
 For the second case, we obtain the same conclusion using  Lemmas 5.7. and 5.9
and checking all the possible contractions.

 Case (2) : $E \notin \bsa$.
 Case by case checking (based on  Remark 4.4 and applications of
 the lemmas 5.6 through 5.10) shows that only case (A) occurs and that
 $(K _{S_1} + \bda _ 1 ) \cdot \hat E = 0$.
\qed
\enddemo

\medpagebreak

The following theorems are straightforward consequences of the above results.
 All the statements are obtained listing all the possible contractions for
the curve not in $\bsa_{X/T}$ or not at the ``feet" and using the lemmas to
show
that only the one listed can occur.

\proclaim{ Theorem 4.5}
 The graphs listed in  Tables 2, 3, 4  are the only possible on the minimal
resolution
 of $\bar S$.
In particular only Hirzebruch-Jung  strings can occur.
The intersecting non log extremal curves need not to have simple normal
crossings,
 but the sum of their coefficient has to equal the sum of the coefficient of
the curves represented by empty circles in the Tables.

\endproclaim
\demo{\bf Proof}
\noindent Set $\phi _2 \cdot \phi _1  = \phi' _2 \cdot \phi' _1 $,  where
$\phi' _1 $ denotes the contractions of $\hat E$ to a smooth point.   Then the
algorithm in the second part of Proposition 2.5. still applies (see note after
4.1).
\qed
\enddemo

\proclaim{\bf Theorem 4.6}
$$ ( K_S + \bda) \cdot \bda =   [ ( K _ {\bar S  } + \bar \bda  )
      \cdot \bar \bda  - \sum \Sb p \in U \endSb (1 - 1/ {n_p}) ] \text{  ,}$$
 where $ U = \{ p / p \text{ is a singular
points of } \bar S \} $ and the exceptional graph of $p$ is a Hirzebruch-Jung
string
of type $ <n_p , q_p>$
\endproclaim

\demo{\it Proof}
{}From Proposition 4.1. we know that we can re-factor
the contraction morphism, as decribed in \S 2.1. and that Proposition 2.2.
still applies. \qed
\enddemo

\proclaim{  Proposition [G1]}
Let $\pi :X \rightarrow S$ be an elliptic fibration between smooth varieties
and $\bsa$ a divisor of simple normal crossings.  Let $n$ be the dimension of
$X$.  Then:
\roster
\item $ \chi (\Cal O  _X ) =  \chi ( S, \ \Cal O _S (- \pi _* K_ {X / S}) ) +
\chi (S, \Cal O_S) $
\item  If $n=3$ then, $\chi (\Cal O _X) = - \frac { (K_S + \pi _* K_ {X / S} )
\cdot \pi _* K_ {X / S} } {2} \ .$
\endroster
\endproclaim
Combining the above Proposition with Theorem 4.6 we have:

\proclaim{ Corollary 4.7} Assume that all the exceptional curves for $\phi$
which are image of multiple fiber divisors have zero intersection with  curves
image of  non multiple singular divisor.
Then:
$$ \chi ( \Cal O _X) = - \frac{1}{2} \{ [ K _ {\bar S  } + \phi _* \pi _*
(K_{X/S}) ]
      \cdot \phi_* \pi _*(K_{X/S}) - \sum \Sb p \in V \endSb (1 - 1/ {n_p})
\},$$
 where $ V = \{ p \in U \text{ such that the fiber over the general point of  }
\phi ^{-1} (p) \text{ is not multiple } \}$.

\endproclaim

\bigpagebreak

\heading{\bf  \S 5 The Lemmas}
\endheading

\noindent{\bf  Preparation for the first step} (Contractions to smooth points.)

The notation is as in (0.11): $\{ E_1, \dots,  E _{k,k-1} \dots E _{r,r-1} \}$
 is a string of log extremal curves on a smooth surface $S$ and $\psi _i :
S_{i - 1} \rightarrow S_i$ , $ 1 \leq i \leq r$
the contraction of the log extremal curve $E_i \subset S _{i-1}$.

\smallpagebreak

\proclaim{\bf Proposition 5.1}
$$\aligned
K_S \cdot E_{k, k-1} &= \psi _1 ^* (K_1) \cdot E_{k, k-1} = K_1 \cdot E_{k,
k-2} = \cdots = K_{k-1} E_k + \beta _{k-1} = \beta _k E_k ^2 + \beta _{k-1}\\
\bda \cdot E_{k, k-1} &= \psi _1 ^* (\bda_1) \cdot E_{k, k-1} = \bda_1 \cdot
E_{k, k-2} =  \cdots = \bda_{k-1} E_k + \alpha _{k-1} = \alpha _k E_k ^2 +
\alpha _{k-1}\\
 -b _k &= ( E  _{k,k-1} ) ^2 = \dots  =( E   _{k,2} ) ^2 = ( E  _{k,1} ) ^2
\endaligned$$
\endproclaim
\demo {\it Proof}
It follows directly from the definitions.
\enddemo

 \proclaim{ \bf Lemma 5.2}
\roster
\item $(E_ k) ^ 2 = - f _k / g _k$
\item $\psi ^* _ {k-1} ( E_k) = E _{k,1} + ( g_{k-1} / f_{k-1} )  E _{k-1}$.
\endroster
\endproclaim
\demo{\it Proof}
$\psi ^* _{k-1} ( E_k) = E  _{k,1} + \mu _k E _{k-1}$, for some rational number
$\mu _k$.
 Intersecting both sides of the equation  first with $E_ {k-1}$ and then with
$E_{k,1}$,
we get $  \mu _k = - 1 / ( E _{k-1}) ^2$ and
 $(E_k) ^2 = - b_k - 1 / (E _{k-1}) ^2 . $

Since $(E _1) ^2 = -  b _1 = - f _1 / g _1 $, (1) holds for $k =1$.
 Assuming by induction $(E_ {k-1}) ^ 2 = - f _{k-1} / g _{k-1}$ we have $(E_k)
^2 = - \frac{f_k}{g_k}$ (0.4)
and $\psi ^* _ {k-1} ( E_k) = E _{k,1} + \frac{ g_{k-1}}{ f_{k-1} }  E _{k-1}$.
 \qed

\enddemo

\proclaim{ \bf Lemma 5.3 }
\roster
\item $\beta _k = ( 2 - b_k + \beta _{k-1} ) \ \frac{g_k}{f_k} , \ \forall k
\geq 2$.
\item $\beta _k = ( g _k  + 1 - f _k) / f _k \forall k$.
\endroster
\endproclaim

\demo{\it Proof}
(1) Since all the curves in the string are rational,

$ -2 + b_k = K_S \cdot E _{k ,k-1 } \overset (0.9) \to = \beta _k  ( E _k )^2
+ \beta _{k - 1} \ . $

Thus $ \beta _k = \frac{2 - b_k + \beta _{k-1}}{E_k ^2}= (2 - b_k + \beta
_{k-1}) \  \frac{ g_k}{f_k}$.

\smallpagebreak

It follows also that $\beta  _1 = ( 2 - b _1 ) / b _1 \overset  (0.4) \to = (g
_1 + 1 - f _1 ) / f _1 $.

 Thus (2) holds for $k =1$.

(2) Assume by induction    $\beta _{k-1} = ( g _{k-1}  + 1 - f _{k-1} ) / f
_{k-1} $ and
prove that the formula holds for $\beta _k$.
Using part (1), we have:
$$\aligned \beta _k & = [ 2 - b_k + ( g _{k-1} + 1 - f _ {k-1} ) / f_{k-1} ]
(g_k / f_k)\\
   \beta _k & \overset( 0.10) \to =  ( 2 g _k - b _k f _{k-1} + g _{k-1} + 1 -
g_k ) / f_ k  \\
    \beta _k & =  ( g _k + 1 - f _k ) / f _k  \ \qed
\endaligned
$$
\enddemo

\vskip 0.3in

\noindent {\bf Preparation for the second step} (Contractions of
Hirzebruch-Jung strings.)

\medpagebreak

Now we consider a string of log extremal Hirzebruch-Jung curves  $\{ E_1 ,
E_{2,1} , \dots , E_{r, r-1} \}$ (as in 0.9).

\smallpagebreak

 \noindent {\bf Remark 5.4: }

\roster
\item  Set  $\bda \cdot  E  _{k,k-1}= m_k$.  Then  $m_k$ is an integer (0.14,
part 3), and
$$ m_ k = \bda  _ 1 \cdot E  _ {k,k-2} = \dots =
\bda _ {k-2} \cdot E  _{k,2} =
\bda _ {k-1} \cdot E  _{k,1}
        + \alpha _ {k-1} = \alpha_k \cdot (E _k) ^2 + \alpha _{k-1}.$$
\item  $\alpha _k  +       \beta_k  > 0 , \ \forall k$.
\item Consequently, $0 < \alpha <1 $, since $ \ \forall k \ \beta _k \leq 0, $
and  $\alpha _k  < 1$ (2.1.).
\endroster

The goal of the first part of this session is to prove Proposition 2.2,
which is part of the combinatorial data needed for the contraction algorithm.

We need the following lemma.

\proclaim{\bf  Lemma 5.5}
\roster
\item
$\alpha _ k = 1 - g _k / f _ k  \text{ , } \forall k$.
\item  $\alpha _ k + \beta _ k = 1 /f _ k$.
\endroster
\endproclaim

\demo{\it Proof}  (1) The crucial point in the proof of the Lemma is that $m_k$
is an integer (0.14, part 3).

 We will show that the Lemma holds for $k=1$ and then use induction on $k$.
 Set $\alpha_ 1 =\frac{p_1}{q_1},$  with  $p_1, q_1$ relatively prime  and $$ 0
< p_1 < q_1  \ (5.4) ;  $$ then $m_1= - \alpha _1  b _1 =  - p_1
\frac{b_1}{q_1}$ (note that $ m_1  < 0$) .

 Since $m_1$ must be an integer, $ b_1 = l q _1$, for some positive integer
$l$; hence
    $$\alpha _1 + \beta _1 = \alpha _1 + \frac{2 - b_1}{b_1} > 0 \ \text{ and }
p _1 l + 2 - l q _1 > 0 .$$
It follows that $ \ l q _1 - 2 < l p _1 < l q _1$ and thus
      $ l p _ 1 = l q _1  - 1 $.

 Therefore :
$$\left\{
\aligned l& = 1\\
         p _ 1 & = q _ 1 - 1 \endaligned
\right.
\quad
\Rightarrow
\quad
\left\{
\aligned q_1& = b_1\\
         p _ 1 & = b _ 1 - 1 \endaligned
\right.
 $$

Hence $\alpha _1 = ( b_ 1 - 1 ) / b _ 1 = 1 - g _ 1 / f _ 1 $ and
 the lemma holds for $ k = 1 $.

\smallpagebreak

 Assume by induction
$  \alpha _{ k -1}  = 1 - g _ {k-1} / f _ {k-1} $ and write  $\alpha_ k = p_ k
/ q _k$, with $p_k, \ q_k$ relatively prime integers and $ 0 < p_k < q_k$
(5.4).

By 5.4 and 5.2. $m_k =  \alpha _ k (E_k)^2  +  \alpha _  {k-1} = - \alpha _ k
\frac{f _ k }{ g _k } +  \alpha _  {k-1}$.
  Then:
$$\align m _ k& = - (1/ f _{k-1}) [ (p_k / q_k )( b_k f_{k-1} - g _{k-1} ) + g
_{k-1} ] + 1\\
   m _ k& = - (1/ f _{k-1}) [ (p_k / q_k ) f _k + g _{k-1} ] + 1 . \tag 5.5.1
\endalign
 $$
$m_k$ must be an integer, and thus $f_k = l_k  q_k$ for some positive integer $
l_k$;
 5.5.1 becomes then:
 $$m_k = - \frac{1}{f_{k-1}} \  ( p _k l_k  + g _{k-1} ) + 1. \tag 5.5.2$$
By (5.3)
$$\alpha _k + \beta _k >0 \Longleftrightarrow
\alpha _ k f _ k + ( 2- b _k) g_k + ( g_{k-1} + 1 - f _{k-1} ) > 0,  \ \forall
k \geq 2 $$
 that is
$$ p_k l_k + (2- b_k) g_k + (g_{k-1} +1 - f_{k-1}) >0. $$

By (0.10) this last equation becomes
$$p_k l_k > b _k f _{k-1} - f _{k-1} - 1 - g _{k-1};$$
 also $l_ k p_k < f_k \overset (5.4) \to= b_{k-1}f_{k-1} - g_{k-1}$.

Combining the last two inequalities we get:
$$\aligned ( b_k -1)f_{k-1}  - 1  &< p _k l_k + g _{k-1} < b_k f_{k-1}
\quad   \text{ i.e.} \\
 ( b _k  - 1 ) f_{k-1} &\leq p _k l_k + g _{k-1} < b_k f_{k-1}.\endaligned
$$
Since $m_k$ is an integer, $f_{k-1}$ must divide $p _k l_k + g _{k-1}$ (5.5.2);
thus:
$ p_k l_k  = (b_k - 1 ) f _{k-1} - g _{k-1} $, that is
 $$  p_k l_k  = f_k - f_{k-1} \tag 5.5.3$$
Dividing both side of (5.5.3) by $q_k l_k = f _k$, we have:
$$\aligned\alpha _ k& = \frac{p_k l_k}{q_k l_k} = 1 - \frac{f_{k-1} }{ f_k}
\text{ or}\\
   \alpha _ k& = 1\ - \ g _k  / f _k . \endaligned
$$

(2) $\alpha _ k + \beta _ k = 1 /f _ k$ follows from  part (1), 5.3. and
(0.10).  \qed
\enddemo

\bigpagebreak

{\bf Preparation for the third step}

\medpagebreak
First the consider the case of a log extremal curve at the feet of $s$
Hirzebruch-Jung strings $\{ \Cal E_1, \cdots \Cal E_s \}, \ s = 2, 3$ (see
4.3).

\proclaim{ Lemma  5.6}
Let $\{ \Cal E_1,\cdots  ,  \Cal E_s \}$
be $s =2, 3$ disjoint Hirzebruch-Jung  strings, and $F_0$ a smooth rational
curve
at the feet of $ \{ \Cal E _1,\cdots  ,\Cal E_s \}$.

Then the image of $F_0$ is not log extremal on $T$.
\endproclaim

\demo{\it Proof} We will prove the statement for the case $s=2$. The same
argument works for the case $s=3$,
but the notation is more complicated.

Assume that such a $F_0$ exists and derive a contradiction.
The proof is along the same line of 5.5. ($\bda \cdot F_0$ is an
integer and the contraction of the image of $F_0$ is log extremal.)

Let $\psi _1: S \rightarrow T_1$ be the contraction of the string $ \Cal E_ 1$
and $\psi _2 : T_1 \rightarrow T_2$ the contraction of the string $\Cal E _2$.
Set $F_j = {(\psi _j)}_* (F_{j-1})$ and denote by
  $\psi _3  : T_2 \rightarrow T_3$ the contraction of $F_2$; set
 $   K _{T_2} = \psi _3 ^* (K_ {T_3}) + \zeta F_2$ and $ \bda _{T_2} = \psi _
3^* (\bda_ {T_3}) + \eta F_2$.

Then
$$  ( K _S + \bda ) \cdot F_0 = (\eta + \zeta ) ( F _2 ) ^2  +  \frac{1}{f
_{r_1}} + \frac{ 1}{ f _{r_2}}. $$  If we set $(F_0)^2 = -b$, it is easy to see
that
$$
\left\{
\aligned (F_2) ^2 &= -b +\frac{g_{r_1}}{f _{r_1}} + \frac{g_{r_2}}{ f _{r_2}}
\\
 \zeta &=  [b -2 -  (\frac{g_{r_1} - f_{r_1} +1}{f _{r_1}} + \frac{ g_{r_2} -
f_{r_2} +1}{f _{r_2}})]
\ \big/ \ (  \frac{g_{r_1} }{ f_{r_1}} + \frac{g_{r_2} }{ f_{r_2} -b})
\endaligned
\right.
$$
That is $\zeta = -1 + \frac{ f_{r_1} + f_{r_2}  }{ b_1 f_{r_1} f_{r_2} -
f_{r_2} g_{r_1} - f_{r_1} g_{r_2}} .$
Set $\eta = p/q$.
Since $\bda \cdot F_0 \in \Bbb Z $,
$$ \frac{1}{(f_{r_1} f_{r_2}) } \{ \frac{p}{q }(bf_{r_1} f_{r_2}  -g_{r_2}
f_{r_1} -g_{r_1} f_{r_2}) + g_{r_1} f_{r_2} + g_{r_2} f_{r_1} \}$$
must be an integer.  Therefore $ bf_{r_1} f_{r_2} -g_{r_2} f_{r_1} -g_{r_1}
f_{r_2} = A q $, for some integer  $A$; thus $\zeta = \frac{- Aq + f_{r_1} +
f_{r_2}}{Aq}$ and $ \frac{1}{(f_{r_1} f_{r_2)} } \{ \frac{p}A + g_{r_1} f_{r_2}
+ g_{r_2} f_{r_1} \}$ is an integer.

The singularities are log terminal ($\zeta  >-1$), thus $ Aq  $ is a positive
integer.
As in the proof on  5.5 have  $ -\zeta < \eta <1 $ and thus:
$-\zeta qA   < pA < q A$, that is
$$\aligned b f_{r_1} f_{r_2} - (f_{r_1} + f_{r_2}) - f_{r_1} g_{r_2} - f_{r_2}
g_{r_1}  & < p A
< p ( b f_{r_1} f_{r_2} - f_{r_2} g_{r_1} - f_{r_1} g_{r_2}) \\
 b f_{r_1} f_{r_2} - (f_{r_1} + f_{r_2}) &< Ap +  f_{r_2} g_{r_1} + f_{r_1}
g_{r_2} < b f_{r_1} f_{r_2} \endaligned
$$
Since $\bda \cdot F_0$ is an integer, $f_{r_1} f_{r_2}$ must divide $Ap +
f_{r_2} g_{r_1} + f_{r_1} g_{r_2}$,
 and thus  such $p$ and $q$ exist if and
only if $f_{r_1} + f_{r_2} > f_1 f_{r_2}$.
This impossible because $ \Cal E _1$ and $\Cal E_2 $ are Hirzebruch-Jung
strings. \qed
\enddemo

\medpagebreak

Thus a
log extremal curve $E$ can never be at the ``feet" of 2 or more exceptional
curves (see also 4.3).

\medpagebreak

 Here starts a series Lemmas and Corollaries which will be used to prove
 that all the log extremal curves contracted by $\varphi _3$ (see \S2)
have self-intersection $-1$ on the minimal resolution.  Again the
 notation is as in (0.9).

\proclaim{\bf Lemma 5.7}
Let $\{ E_1 ,\dots , E_{r-1,r} \}$ be a Hirzebruch-Jung  string and $ \psi : S
\rightarrow  S_r$  the corresponding log extremal contraction. Let $F_0$ be a
curve on $S$ such that  $E_1 \cdot F_0 = 1$.  Set $F_r = (\psi) _* (F_0) $.
Then:
$$(K_S + \bda) \cdot F _0 =  ( K_{S _r} + \bda _r ) \cdot F_r +  \frac{1}{f_1}
+ \sum _{k =1} ^{r-1} \frac{1 }{ f_k f_{k+1}}.$$
\endproclaim

\demo{\it Proof}
Let $\psi _j : S_{j-1} \rightarrow S_j$ be as in (0.9.) and set $F_j = (\psi
_j) _* (F_{j-1}) $.
The statement follows from
$$(K_S + \bda) \cdot F_0 = K_{S _r} + \bda _r +
(\alpha _1 + \beta _1  + \dots + (\alpha _k + \beta _k) E _k \cdot F_{k-1} +
\dots  (\alpha _r + \beta _r) E _r \cdot F_{r-1} $$
and Proposition 5.2. \qed

\enddemo

\proclaim{\bf Corollary 5.8}
In the hypothesis and notation of Lemma 5.7., let $r=2, 3$.
Then $(K_S + \bda ) \cdot F_0 = 0$.

\endproclaim
\demo{\it Proof}
By 5.7, $$0 > (K_S + \bda) \cdot F _0 =  ( K_{S _r} + \bda _r ) \cdot F_r +
\frac{1}{f_1} + \sum _k ^{r-1} \frac{1 }{ f_k f_{k+1}}.$$
\roster
\item r = 2. It's enough to observe :
$$1/f_1 + 1/(f_1 f_2 ) =  1/b_1 + 1/ b_1( b_1 b_2 -1) = b_2 /(b_1 b_2 - 1) < 1
.$$
(In fact,  $ b_2  < (b_1 b_2 - 1) $  unless $b_1 = 1$;  note also that  $(K_S +
\bda) \cdot F _0 $ is a positive integer.)
\item  r = 4.  As above  we have :
$$1/f_1 + 1/(f_1 f_2 ) + 1/ (f_2 f_3 ) = [ b_2 b_3 -1] / [ b_3 ( b_1 b_2 -1) -
b_1 ] < 1 $$
\endroster
\enddemo

\proclaim{\bf Lemma 5.9}
Let $\{ E_1 \text { , } E_{1,2} \text { , } E_{2,3} \}$ be a Hirzebruch-Jung
string and $ \psi : S \rightarrow  S_3$  the corresponding log extremal
contraction. Let $F _0$ be a curve on $S$ such that  $E_{1,2} \cdot F_0 = 1$
and set: $\psi _* (F_0) =  F_3$.  Then:
$$(K_S + \bda) \cdot F_0 =  ( K_{ S_3 }+ \bda _3) \cdot  F _3 + \frac{1}{f_2} +
\frac{1}{f_3}  + \frac{f_1}{f_2}$$

\endproclaim

\demo{\it Proof}
 Procede as in the proof of Lemma 5.7. \qed

\enddemo

\medpagebreak

The following proposition shows that the strict transform on $S_1$ of any log
extremal curve contracted by
 $\varphi _3$ has intersection number $0$ with $K _{S_1} + \bda _1$.

\proclaim{\bf Lemma 5.10}
Let $\{ E_1 \dots E_{r-1,r} \}$ be as in 2.2.
Assume that $E$ is a log extremal curve on $S_r$ and denote by  $\hat E$ the
strict transform of $E$.
Then if $\hat E \cdot E_{r-1, r} =1$, we have:
 $$(K _S + \bda) \cdot \hat E = 0. $$
\endproclaim

\demo{\it Proof}
Without loss of generality (2.2.), we can assume that $\hat E^2 = -1$ and
 $(K_S + \bda ) \cdot \hat E \geq 0$.
Since $$(K_S + \bda ) \cdot \hat E = (K_{S_r} + \bda _r )  E +\alpha _ r +
\beta _ r  $$
we have:
$$0 >  (K_{S_r} + \bda _r )  E = (K_S + \bda ) \cdot \hat  E - 1/f _ r \text{ ;
} $$
Now $(K_S + \bda ) \cdot E$ has to be zero, because $ 1/ f_r < 1 $ and $(K_S +
\bda ) \cdot E \in \Bbb Z $.\qed
\enddemo

 \enddocument

\documentstyle[11pt]{article}

\oddsidemargin=0.5in
\evensidemargin=0in
\textwidth=9in

\headheight=0pt
\headsep=0pt
\topmargin=1in
\textheight=9in

\begin{document}

{\large}

\vskip 0.5in

\setlength{\unitlength}{.8in}

r=1

\begin{picture}(7,2)(0,0)

\thicklines

\put(0,1){\circle*{.075}}

\put(0.3,1){\makebox(.25,0){ (-3)}}
\put(0,.5){\makebox(.25,.25){$IV^*$}}

\thicklines

\put(2,1){\circle*{.075}}

\put(2.3,1){\makebox(.25,0){ (-6)}}
\put(2,.5){\makebox(.25,.25){$II^*$}}

\thicklines

\put(4,1){\circle*{.075}}

\put(4.3,1){\makebox(.25,0){(-4)}}

\put(4,.5){\makebox(.25,.25){$III^*$ }}

\thicklines

\put(6,1){\circle*{.075}}

\put(6.3,1){\makebox(.25,0){(-2) }}

\put(6.3,.5){\makebox(.25,.25){$I_0^ *$}}

\put(6.35,1.4){\makebox(.25,.25){$I_a$}}

\end{picture}

\bigskip

\begin{picture}(7,2)(0,0)

\thicklines
\put(0,1){\circle*{.075}}

\put(0,1){\line(1,0){.4625}}

\put(0.5,1){\circle{.075}}

\put(0.7,1){\makebox(.25,0){ (-5)}}

\put(0,.5){\makebox(.25,.25){$II^*$}}

\put(0.4,.5){\makebox(.25,.25){$II$}}

\thicklines
\put(2,1){\circle*{.075}}

\put(2,1){\line(1,0){.4625}}

\put(2.5,1){\circle{.075}}

\put(2.7,1){\makebox(.25,0){ (-4)}}

\put(2,.5){\makebox(.25,.25){$II^*$}}

\put(2.4,.5){\makebox(.25,.25){$IV$}}

\thicklines
\put(4,1){\circle*{.075}}

\put(4,1){\line(0,1){.4625}}

\put(4,1.5){\circle{.075}}

\put(4,1){\line(1,0){.4625}}

\put(4.5,1){\circle{.075}}

\put(4.7,1){\makebox(.25,0){ (-4)}}

\put(4,.5){\makebox(.25,.25){$II^*$}}

\put(4.04,1.5){\makebox(.25,.25){$II$}}

\put(4.4,.5){\makebox(.25,.25){$II$}}

\thicklines

\put(6,1){\circle*{.075}}

\put(6,1){\line(1,0){.4625}}

\put(6.5,1){\circle{.075}}

\put(6.7,1){\makebox(.25,0){ (-3)}}

\put(6,.5){\makebox(.25,.25){$II^*$}}

\put(6.4,.5){\makebox(.25,.25){$I_0 ^*$}}

\end{picture}

\bigskip

\begin{picture}(7,2)(0,0)

\thicklines
\put(0,1){\circle*{.075}}

\put(0,1){\line(0,1){.4625}}

\put(0,1.5){\circle{.075}}

\put(0,1){\line(1,0){.4625}}

\put(0.5,1){\circle{.075}}

\put(0.7,1){\makebox(.25,0){ (-3)}}

\put(0,.5){\makebox(.25,.25){$II^*$}}

\put(0.04,1.5){\makebox(.25,.25){$IV$}}

\put(0.4,.5){\makebox(.25,.25){$II$}}

\thicklines
\put(1.5,1){\circle{.075}}

\put(1.5,1){\line(1,0){.4625}}

\put(2,1){\circle*{.075}}

\put(2,1){\line(0,1){.4625}}

\put(2,1.5){\circle{.075}}

\put(2,1){\line(1,0){.4625}}

\put(2.5,1){\circle{.075}}

\put(2.7,1){\makebox(.25,0){ (-3)}}

\put(1.5,.5){\makebox(.25,.25){$II$}}

\put(2,.5){\makebox(.25,.25){$II^*$}}

\put(2.04,1.5){\makebox(.25,.25){$II$}}

\put(2.4,.5){\makebox(.25,.25){$II$}}

\thicklines
\put(4,1){\circle*{.075}}

\put(4,1){\line(1,0){.4625}}

\put(4.5,1){\circle{.075}}

\put(4.7,1){\makebox(.25,0){ (-3)}}

\put(4,.5){\makebox(.25,.25){$III^*$}}

\put(4.5,.5){\makebox(.25,.25){$III$}}

\thicklines
\put(6,1){\circle*{.075}}

\put(6,1){\line(0,1){.4625}}

\put(6,1.5){\circle{.075}}

\put(6,1){\line(1,0){.4625}}

\put(6.5,1){\circle{.075}}

\put(6.7,1){\makebox(.25,0){ (-2)}}

\put(6,.5){\makebox(.25,.25){$III^*$}}

\put(6.5,.5){\makebox(.25,.25){$III$}}

\put(6.04,1.5){\makebox(.25,.25){$III$}}

\end{picture}

\bigskip

\begin{picture}(7,2)(0,0)

\thicklines
\put(0,1){\circle*{.075}}

\put(0,1){\line(0,1){.4625}}

\put(0,1.5){\circle{.075}}

\put(0,1){\line(1,0){.4625}}

\put(0.5,1){\circle{.075}}

\put(0.7,1){\makebox(.25,0){ (-2)}}

\put(0,.5){\makebox(.25,.25){$II^*$}}

\put(0.5,.5){\makebox(.25,.25){$I_0 ^*$}}

\put(0.04,1.5){\makebox(.25,.25){$II$}}

\thicklines
\put(2,1){\circle*{.075}}

\put(2,1){\line(1,0){.4625}}

\put(2.5,1){\circle{.075}}

\put(2.7,1){\makebox(.25,0){ (-2)}}

\put(2,.5){\makebox(.25,.25){$III^*$}}

\put(2.5,.5){\makebox(.25,.25){$I_0 ^*$}}

\thicklines

\put(3.5,1){\circle{.075}}

\put(3.5,1){\line(1,0){.4625}}

\put(4,1){\circle*{.075}}

\put(4,1){\line(0,1){.4625}}

\put(4,1.5){\circle{.075}}

\put(4,1){\line(1,0){.4625}}

\put(4.5,1){\circle{.075}}

\put(4.7,1){\makebox(.25,0){ (-2)}}

\put(3.5,.5){\makebox(.25,.25){$II$}}

\put(4,.5){\makebox(.25,.25){$II^*$}}

\put(4.5,.5){\makebox(.25,.25){$II$}}

\put(4.04,1.5){\makebox(.25,.25){$IV$}}

\thicklines

\put(5.5,1){\circle{.075}}

\put(5.5,1){\line(1,0){.4625}}

\put(6,1){\circle*{.075}}

\put(6,1){\line(1,1){.4625}}

\put(6,1){\line(-1,1){.4625}}

\put(5.5,1.5){\circle{.075}}

\put(6.5,1.5){\circle{.075}}

\put(6,1){\line(1,0){.4625}}

\put(6.5,1){\circle{.075}}

\put(6.7,1){\makebox(.25,0){ (-2)}}

\put(5.5,.5){\makebox(.25,.25){$II$}}

\put(6,.5){\makebox(.25,.25){$II^*$}}

\put(6.5,.5){\makebox(.25,.25){$II$}}

\put(6.55,1.4){\makebox(.25,.25){$II$}}

\put(5.55,1.4){\makebox(.25,.25){$II$}}

\end{picture}

\newpage

\vskip 0.5in

\begin{picture}(7,2)(0,0)

\thicklines
\put(0,1){\circle*{.075}}

\put(0,1){\line(1,0){.4625}}

\put(0.5,1){\circle{.075}}

\put(0.7,1){\makebox(.25,0){ (-2)}}

\put(0,.5){\makebox(.25,.25){$II^*$}}

\put(0.4,.5){\makebox(.25,.25){$IV^*$}}

\thicklines
\put(2,1){\circle*{.075}}

\put(2,1){\line(1,0){.4625}}

\put(2.5,1){\circle{.075}}

\put(2.7,1){\makebox(.25,0){ (-2)}}

\put(2,.5){\makebox(.25,.25){$IV^*$}}

\put(2.4,.5){\makebox(.25,.25){$IV$}}

\thicklines
\put(4,1){\circle*{.075}}

\put(4,1){\line(0,1){.4625}}

\put(4,1.5){\circle{.075}}

\put(4,1){\line(1,0){.4625}}

\put(4.5,1){\circle{.075}}

\put(4.7,1){\makebox(.25,0){ (-2)}}

\put(4,.5){\makebox(.25,.25){$II^*$}}

\put(4.04,1.5){\makebox(.25,.25){$IV$}}

\put(4.4,.5){\makebox(.25,.25){$IV$}}

\thicklines

\put(6,1){\circle*{.075}}

\put(6,1){\line(0,1){.4625}}

\put(6,1.5){\circle{.075}}

\put(6,1){\line(1,0){.4625}}

\put(6.5,1){\circle{.075}}

\put(6.7,1){\makebox(.25,0){ (-2)}}

\put(6,.5){\makebox(.25,.25){$IV^*$}}

\put(6.04,1.5){\makebox(.25,.25){$II$}}

\put(6.4,.5){\makebox(.25,.25){$II$}}

\end{picture}

\bigskip

r = 2.

\begin{picture}(7,2)(0,0)

\thicklines
\put(0,1){\circle*{.075}}

\put(0.5,1){\line(0,1){.4625}}

\put(0.5,1.5){\circle{.075}}

\put(0,1){\line(1,0){.4625}}

\put(0.5,1){\circle*{.075}}

\put(0.7,1){\makebox(.25,0){ (-2,-3)}}

\put(0,.5){\makebox(.25,.25){$II^*$}}

\put(0.06,1.4){\makebox(.25,.25){$II$}}

\put(0.4,.5){\makebox(.25,.25){$IV*$}}

\thicklines
\put(2,1){\circle*{.075}}

\put(2.5,1){\line(0,1){.4625}}

\put(2.5,1.5){\circle{.075}}

\put(2,1){\line(1,0){.4625}}

\put(2.5,1){\circle*{.075}}

\put(2.7,1){\makebox(.25,0){ (-3,-2)}}

\put(2,.5){\makebox(.25,.25){$II^*$}}

\put(2.06,1.4){\makebox(.25,.25){$II$}}

\put(2.4,.5){\makebox(.25,.25){$I_0 *$}}

\thicklines
\put(4,1){\circle*{.075}}

\put(4.5,1){\line(0,1){.4625}}

\put(4.5,1.5){\circle{.075}}

\put(4,1){\line(1,0){.4625}}

\put(4.5,1){\circle*{.075}}

\put(4.5,1){\line(1,0){.4625}}

\put(5,1){\circle{.075}}

\put(5.2,1){\makebox(.25,0){ (-3,-2)}}

\put(4,.5){\makebox(.25,.25){$II^*$}}

\put(4.56,1.4){\makebox(.25,.25){$I_1$}}

\put(4.4,.5){\makebox(.25,.25){$I_0  ^*$}}

\put(5,.5){\makebox(.25,.25){$I_1$}}

\end{picture}

\bigskip

(-2,-2)

\begin{picture}(7,2)(0,0)

\thicklines
\put(0,1){\circle*{.075}}

\put(0.5,1){\line(0,1){.4625}}

\put(0.5,1.5){\circle{.075}}

\put(0,1){\line(1,0){.4625}}

\put(0.5,1){\circle*{.075}}

\put(0,.5){\makebox(.25,.25){$II^*$}}

\put(0.06,1.4){\makebox(.25,.25){$I_0 ^*$}}

\put(0.4,.5){\makebox(.25,.25){$IV*$}}

\thicklines
\put(1.5,1){\circle*{.075}}

\put(2,1){\line(0,1){.4625}}

\put(2,1.5){\circle{.075}}

\put(1.5,1){\line(1,0){.4625}}

\put(2,1){\circle*{.075}}

\put(1.5,.5){\makebox(.25,.25){$III^*$}}

\put(2.1,1.4){\makebox(.25,.25){$III$}}

\put(1.9,.5){\makebox(.25,.25){$I_0 ^*$}}

\thicklines
\put(3,1){\circle*{.075}}

\put(3.5,1){\line(0,1){.4625}}

\put(3.5,1.5){\circle{.075}}

\put(3,1){\line(1,0){.4625}}

\put(3.5,1){\circle*{.075}}

\put(3.5,1){\line(1,0){.4625}}

\put(4,1){\circle{.075}}

\put(3,.5){\makebox(.25,.25){$II^*$}}

\put(3.56,1.4){\makebox(.25,.25){$IV$}}

\put(3.4,.5){\makebox(.25,.25){$IV^ *$}}

\put(4,.5){\makebox(.25,.25){$II$}}

\thicklines
\put(5.5,1){\circle*{.075}}

\put(5.5,1){\line(1,0){.4625}}

\put(6,1){\circle*{.075}}

\put(5.5,.5){\makebox(.25,.25){$IV^*$}}

\put(6,.5){\makebox(.25,.25){$IV$}}

\end{picture}

\bigskip

\begin{picture}(7,2)(0,0)

\thicklines
\put(0,1){\circle*{.075}}

\put(0.5,1){\line(0,1){.4625}}

\put(0.5,1){\line(1,1){.4625}}

\put(0.5,1.5){\circle{.075}}

\put(1,1.5){\circle{.075}}

\put(0,1){\line(1,0){.4625}}

\put(0.5,1){\circle*{.075}}

\put(0.5,1){\line(1,0){.4625}}

\put(1,1){\circle{.075}}

\put(0,.5){\makebox(.25,.25){$II^*$}}

\put(0.56,1.4){\makebox(.25,.25){$II$}}

\put(0.4,.5){\makebox(.25,.25){$IV ^ *$}}

\put(1,.5){\makebox(.25,.25){$II$}}

\put(2,1){\circle*{.075}}

\put(2,1){\line(0,1){.4625}}

\put(2,1.5){\circle{.075}}

\put(2.5,1){\line(0,1){.4625}}

\put(2.5,1.5){\circle{.075}}

\put(2,1){\line(1,0){.4625}}

\put(2.5,1){\circle*{.075}}

\put(2,.5){\makebox(.25,.25){$II^*$}}

\put(2.06,1.4){\makebox(.25,.25){$II$}}

\put(2.4,.5){\makebox(.25,.25){$I_0 ^ *$}}

\put(2.56,1.4){\makebox(.25,.25){$II$}}

\thicklines
\put(3.5,1){\circle*{.075}}

\put(3.5,1){\line(0,1){.4625}}

\put(3.5,1.5){\circle{.075}}

\put(4,1){\line(0,1){.4625}}

\put(4,1.5){\circle{.075}}

\put(3.5,1){\line(1,0){.4625}}

\put(4,1){\line(1,0){.4625}}

\put(4,1){\circle*{.075}}

\put(4.5,1){\circle{.075}}

\put(3.5,.5){\makebox(.25,.25){$II^*$}}

\put(3.56,1.4){\makebox(.25,.25){$II$}}

\put(3.9,.5){\makebox(.25,.25){$I_0 ^ *$}}

\put(6,.5){\makebox(.25,.25){$\sum =1/4$}}

\put(4.5,.5){\makebox(.25,.25){$II^*$}}

\thicklines

\put(5,1){\circle*{.075}}

\put(5,1){\line(0,1){.4625}}

\put(5,1.5){\circle{.075}}

\put(5.5,1.5){\circle{.075}}

\put(5,1){\line(1,0){.4625}}

\put(5.5,1){\line(0,1){.4625}}

\put(5.5,1){\line(1,0){.4625}}

\put(5.5,1){\circle*{.075}}

\put(6,1){\circle{.075}}

\put(5,.5){\makebox(.25,.25){$III^*$}}

\put(5.4,.5){\makebox(.25,.25){$I_0 ^ *$}}

\put(6,.5){\makebox(.25,.25){$\sum =1/4$}}

\end{picture}

\newpage

\vskip 0.5in

All (-2) strings.

\bigskip

r = 3

\begin{picture}(7,2)(0,0)

\thicklines
\put(0,1){\circle*{.075}}

\put(0,1){\line(1,0){.4625}}

\put(0.5,1){\circle*{.075}}

\put(0.5,1){\line(1,0){.4625}}

\put(1,1){\circle*{.075}}

\put(1,1){\line(0,1){.4625}}

\put(1,1.5){\circle{.075}}

\put(0,.5){\makebox(.25,.25){$II^*$}}

\put(0.4,.5){\makebox(.25,.25){$IV*$}}

\put(0.9,.5){\makebox(.25,.25){$I_0 ^*$}}

\put(1.06,1.4){\makebox(.25,.25){$IV$}}

\thicklines
\put(1.5,1){\circle*{.075}}

\put(1.5,1){\line(1,0){.4625}}

\put(2,1){\circle*{.075}}

\put(2,1){\line(1,0){.4625}}

\put(2.5,1){\circle*{.075}}

\put(2.5,1){\line(0,1){.4625}}

\put(2.5,1.5){\circle{.075}}

\put(2.5,1){\line(1,1){.4625}}

\put(3,1.5){\circle{.075}}

\put(1.5,.5){\makebox(.25,.25){$II^*$}}

\put(1.9,.5){\makebox(.25,.25){$IV*$}}

\put(2.4,.5){\makebox(.25,.25){$I_0 ^*$}}

\put(2.56,1.7){\makebox(.25,.25){$\sum = 1/3 $}}

\thicklines
\put(3.5,1){\circle*{.075}}

\put(3.5,1){\line(1,0){.4625}}

\put(4,1){\circle*{.075}}

\put(4,1){\line(1,0){.4625}}

\put(4.5,1){\circle*{.075}}

\put(4.5,1){\line(0,1){.4625}}

\put(4.5,1.5){\circle{.075}}

\put(4.5,1){\line(1,1){.4625}}

\put(5,1.5){\circle{.075}}

\put(4.5,1){\line(1,0){.4625}}

\put(5,1){\circle{.075}}

\put(3.5,.5){\makebox(.25,.25){$II^*$}}

\put(3.9,.5){\makebox(.25,.25){$IV*$}}

\put(4.4,.5){\makebox(.25,.25){$I_0 ^*$}}

\put(4.56,1.4){\makebox(.25,.25){$I_1$}}

\put(4.9,.5){\makebox(.25,.25){$II$}}

\end{picture}

\begin{picture}(7,2)(0,0)

\thicklines
\put(0,1){\circle*{.075}}

\put(0,1){\line(1,0){.4625}}

\put(0.5,1){\circle*{.075}}

\put(0.5,1){\line(1,0){.4625}}

\put(1,1){\circle*{.075}}

\put(1,1){\line(0,1){.4625}}

\put(1,1.5){\circle{.075}}

\put(1,1){\line(1,0){.4625}}

\put(1.5,1){\circle{.075}}

\put(0,.5){\makebox(.25,.25){$II^*$}}

\put(0.4,.5){\makebox(.25,.25){$IV*$}}

\put(0.9,.5){\makebox(.25,.25){$I_0 ^*$}}

\put(1.06,1.4){\makebox(.25,.25){$I_1$}}

\put(1.4,.5){\makebox(.25,.25){$III$}}

\thicklines
\put(2.5,1){\circle*{.075}}

\put(2.5,1){\line(1,0){.4625}}

\put(3,1){\circle*{.075}}

\put(3,1){\line(1,0){.4625}}

\put(3.5,1){\circle*{.075}}

\put(3.5,1){\line(0,1){.4625}}

\put(3.5,1.5){\circle{.075}}

\put(3.5,1){\line(1,0){.4625}}

\put(4,1){\circle{.075}}

\put(2.5,.5){\makebox(.25,.25){$II^*$}}

\put(3,.5){\makebox(.25,.25){$IV*$}}

\put(3.4,.5){\makebox(.25,.25){$I_0 ^*$}}

\put(3.06,1.4){\makebox(.25,.25){$II $}}

\put(3.9,.5){\makebox(.25,.25){$II$}}

\thicklines

\put(4.5,1){\circle*{.075}}

\put(4.5,1){\line(1,0){.4625}}

\put(5,1){\circle*{.075}}

\put(5,1){\line(0,1){.4625}}

\put(5,1.5){\circle{.075}}

\put(5,1){\line(1,0){.4625}}

\put(5.5,1){\circle*{.075}}

\put(4.5,.5){\makebox(.25,.25){$II ^*$}}

\put(5.06,1.4){\makebox(.25,.25){$II$}}

\put(4.9,.5){\makebox(.25,.25){$IV ^*$}}

\put(5.4,.5){\makebox(.25,.25){$IV $}}

\end{picture}

\bigskip

\begin{picture}(7,2)(0,0)

\thicklines

\put(0,1){\circle*{.075}}

\put(0,1){\line(1,0){.4625}}

\put(0.5,1){\circle*{.075}}

\put(0.5,1){\line(1,0){.4625}}

\put(1,1){\circle*{.075}}

\put(0,.5){\makebox(.25,.25){$III ^*$}}

\put(0.4,.5){\makebox(.25,.25){$I_0 ^* $}}

\put(0.9,.5){\makebox(.25,.25){$III $}}

\end{picture}

\bigskip

r = 4 , r = 5

\begin{picture}(7,2)(0,0)

\thicklines

\put(0,1){\circle*{.075}}

\put(0,1){\line(1,0){.4625}}

\put(0.5,1){\circle*{.075}}

\put(0.5,1){\line(1,0){.4625}}

\put(1,1){\circle*{.075}}

\put(1,1){\line(1,0){.4625}}

\put(1.5,1){\circle*{.075}}

\put(1.5,1){\line(0,1){.4625}}

\put(1.5,1.5){\circle{.075}}

\put(0,.5){\makebox(.25,.25){$II ^*$}}

\put(0.4,.5){\makebox(.25,.25){$IV ^* $}}

\put(0.9,.5){\makebox(.25,.25){$I _0 ^* $}}

\put(1.4,.5){\makebox(.25,.25){$IV$}}

\put(1.06,1.4){\makebox(.25,.25){$II$}}

\thicklines

\put(2,1){\circle*{.075}}

\put(2,1){\line(1,0){.4625}}

\put(2.5,1){\circle*{.075}}

\put(2.5,1){\line(1,0){.4625}}

\put(3,1){\circle*{.075}}

\put(3,1){\line(1,0){.4625}}

\put(3.5,1){\circle*{.075}}

\put(3.5,1){\line(1,0){.4625}}

\put(4,1){\circle*{.075}}

\put(2,.5){\makebox(.25,.25){$II ^*$}}

\put(2.4,.5){\makebox(.25,.25){$IV ^* $}}

\put(2.9,.5){\makebox(.25,.25){$I _0 ^* $}}

\put(3.4,.5){\makebox(.25,.25){$IV$}}

\put(3.9,.5){\makebox(.25,.25){$II$}}

\end{picture}


\vskip 0.4in

 \Refs

\ref\key{\bf  Br}\by  E. Brieskorn \pages 336--358
\paper Rationale singularitaten komplexen Fl\"achen
\yr 1968 \vol 4
\jour Invent. Math.
\endref
\ref\key{\bf G1}\by A. Grassi \pages 287--301
\paper  On minimal model of elliptic threefolds.
\yr 1991 \vol 290
\jour Math. Ann.
\endref
\ref\key{\bf  Ha1}\by  R. Hartshorne  \pages 63--94
\paper Ample vector bundles \vol 29
\yr 1966
\jour Publ. Math. IHES
\endref
\ref\key{\bf Ha2}\by R. Hartshorne
\book Algebraic Geometry
\publ Springer-Verlag \publaddr New York
\yr 1977 \endref
\ref\key{\bf Ka1}\by Y. Kawamata \pages 253--276
\paper Characterization of Abelian varieties.
\yr 1981 \vol 43
\jour Compositio Math.
\endref
\ref\key{\bf Ka2}\by Y. Kawamata \pages 1--24
\paper Kodaira dimension of certain algebraic fiber spaces
\yr 1983 \vol 30
\jour J. Fac. Sci. Tokyo Univ. IA
\endref
\ref\key{\bf Ka3}\by Y. Kawamata
\paper Hodge theory and Kodaira dimension.
\inbook Algebraic Varieties and Analytic Varities
\bookinfo  (S. Iitaka, editor)
Adv. Stud. Pure Math. {\bf 1}
\publ Kinokuniya \publaddr Tokyo \pages 317--327
\yr 1983 \endref
\ref\key{\bf Ka4}\by Y. Kawamata \pages 603--633
\paper The cone of curves of algebraic varieties
\yr 1984 \vol 119
\jour  Ann. of Math.
\endref
\ref\key{\bf Ka5}\by Y. Kawamata \pages 603--633
\paper Pluricanonical systems on minimal algebraic varieties
\yr 1985 \vol 79
\jour  Invent. Math.
\endref
\ref\key{\bf Ka6}\by Y. Kawamata \pages 1--46
\paper Minimal models and the Kodaira dimension of algebraic fiber spaces.
\yr 1985 \vol 363
\jour J. Reine Angew. Math
\endref
\ref\key{\bf Kd}\by K. Kodaira \pages 563--626, 1--40
\paper On compact analytic surfaces II and III
\yr1963\vol77,78
\jour Annals of Math
\endref
\ref\key{\bf Ko}\by J. Koll\'ar \pages 15--36,
\paper Flops
\yr1989\vol113
\jour Nagoya Math. J.
\endref
\ref\key{\bf KMM}\by Y. Kawamata, K. Matsuda, K.Matsuki
\paper  Introduction to the minimal model problem
\inbook Proc. Sym. Alg. Geom. Sendai 1985
\bookinfo Adv. Stud. Pure Math. {\bf 10}
\publ Kinokuniya \publaddr Tokyo \pages 283--360
\yr 1985 \endref
\ref\key{\bf Kz}\by N. Katz \pages 437--443
\paper The regularity theorem in algebraic geometry
\yr 1970 \vol 1
\jour Proc. Internat. Congr. Math., Nice, Gauthier Villars
\endref
\ref\key{\bf Mi}\by Y. Miyaoka \pages 325--332
\paper On the Kodaira dimension of minimal threefolds
\yr 1988 \vol 281
\jour  Math. Ann
\endref
\ref\key{\bf Mo}\by S. Mori \pages 567--588
\paper Flip theorem and the existence of minimal models for 3-folds
\yr 1988 \vol 1
\jour Jour. Amer. Math. Soc.
\endref
\ref\key{\bf Sa}\by  F. Sakai \paper Classification of normal surfaces
\inbook Algebraic Geometry Bowdoin
\bookinfo Proc. Sympos. Pure Math. {\bf 46, 1}
\publ Amer. Math. Soc. \publaddr Providence \pages 451--465
\yr 1987 \endref
\ref\key{\bf Za}\by O. Zariski \pages 560--615
\paper  The theorem of Riemann-Roch for high multiples of an effective divisor
on an algebraic surface
 \yr 1962 \vol 76
\jour Ann. of Math. \endref
\end{document}